\documentclass[a4paper,10pt,oneside,notitlepage]{article}
\usepackage{times}
\usepackage{lmodern}
\usepackage[T1]{fontenc}
\usepackage[utf8]{inputenc}
\usepackage[english]{babel}
\usepackage{amssymb}
\usepackage{amsmath}
\usepackage{mathtools} % it allows to typeset formula in display style without any commands... \begin{dcases}
\usepackage{amsthm}
\usepackage{amsfonts, mathrsfs}
\usepackage{latexsym}
\usepackage[pdftex]{color, graphicx}
\usepackage{makeidx}
\usepackage{sidecap}
\usepackage{verbatim}
\usepackage{geometry}
\usepackage{subfig}
\usepackage{tikz,float}
\usetikzlibrary{shapes,arrows,shadows}
\usepackage[displaymath,mathlines, running]{lineno} % number of line
\usepackage{comment} % label interni
\usepackage[colorinlistoftodos]{todonotes}
\usepackage{esint} % per il simbolo della media integrale
\usepackage{pdfsync}

\newtheorem{thm}{Theorem}[section]
\newtheorem{prop}[thm]{Proposition}
\newtheorem{lemma}[thm]{Lemma}

\newtheorem{remark}[thm]{Remark}
\newtheorem{example}[thm]{Example}

\DeclareMathOperator{\dist}{dist}

\def\R{\mathbb{R}}
\def\Z{\mathbb{Z}}
\def\e{\varepsilon}
\def\cubeTd{Q_{T,d}}
\def\lcubeTd{Q^l_{T,d}}
\def\cubeTk{Q_{T,k}}
\def\loc{{\rm loc}}
\def\hom{{\rm hom}}
\def\weakstar{\overset{\ast}{\rightharpoonup}}

\title{Homogenization of discrete thin structures}

\author{Andrea Braides and Lorenza D'Elia\\ \small Dipartimento di Matematica, Universit\`a di Roma Tor Vergata
\\ \small  via della ricerca scientifica 1, 00133 Roma, Italy}
\date{}                                           % Activate to display a given date or no date

\begin{document}
\maketitle

\abstract{We consider graphs parameterized on a portion $X\subset\mathbb Z^d\times \{1,\ldots, M\}^k$ of a cylindrical subset of the lattice $\mathbb Z^d\times \mathbb Z^k$, and perform a discrete-to-continuum dimension-reduction process for energies defined on $X$ of quadratic type. Our only assumptions are that $X$ be connected as a graph and periodic in the first $d$-directions.  We show that, upon scaling of the domain and of the energies by a small parameter $\e$, the scaled energies converge to a $d$-dimensional limit energy. The main technical points are a dimension-lowering coarse-graining process and a discrete version of the $p$-connectedness approach by Zhikov. 
}

\section{Introduction}
The object of the investigation in this paper is the analysis of discrete thin objects through, at the same time, a discrete-to-continuum and dimension-reduction process. The main focus of our work is the great generality of the geometry of our discrete systems, which we essentially require to be a connected graph periodic in the dimensions that are maintained after a discrete-to-continuum passage. \begin{figure}[h]
\centering
\includegraphics[scale=0.30]{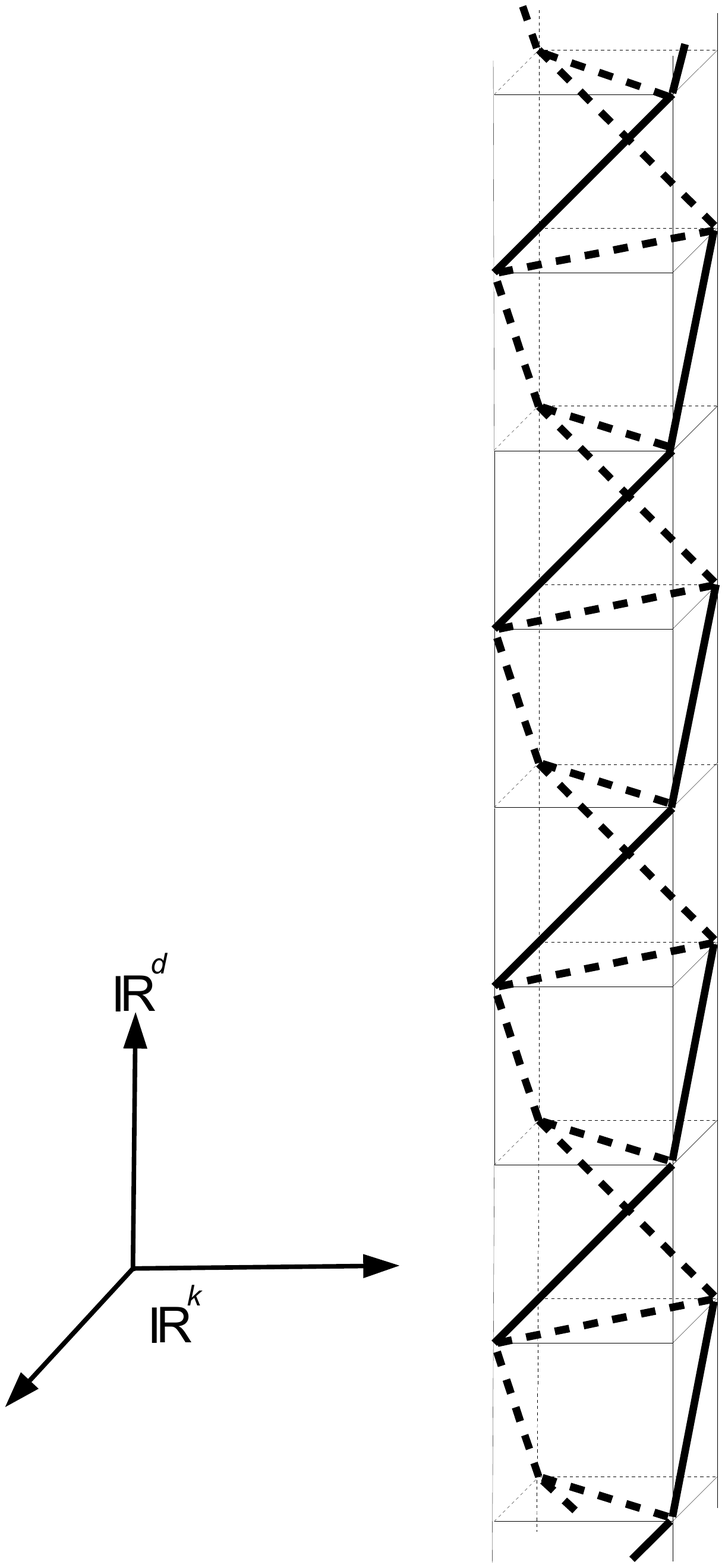}
\caption{A discrete thin object in three dimensions with a one-dimensional behaviour}\label{Fig5larger_vertical}
\end{figure}

An example of the structure that we have in mind is pictured in Fig.~\ref{Fig5larger_vertical}. The thicker black lines (both the solid ones and the ones dashed for graphic purpose) represent connections between nodes of a cubic lattice in $\R^3$.  Equivalently, we may think of the same structure as a network of conducting rods. Note that this object is not trivially a `subgraph' of a function depending of the vertical variable, as it consists of a double helix connected through horizontal bonds.  Nevertheless, it  can be included in a `regular' thin object; in this case, the cylindrical part of $\Z^3$ whose projection on the two-dimensional horizontal plane are the four vertices of a square. Even if no connection is purely vertical, the overall behaviour of such a structure is expected to be that of a vertical one-dimensional object.

With this example in mind, we are going to look at graphs whose nodes are a subset $X$ of $\Z^{d+k}$ periodic of period $T$ in the first $d$ directions (in the example $d=1$, corresponding to the vertical direction, with period $T=2$), bounded in the last $k$ directions (in the example, $k=2$, corresponding to the horizontal directions), so that we may always think that it is contained in $\Z^d\times \{0,\ldots, M-1\}^k$ for some $M\in\mathbb N$. This graph is equipped with a set of edges $\mathcal E\subset X\times X$ which make it connected. This set of edges is supposed to be invariant by the same translations as $X$. 

We are going to show that we may define a continuous $d$-dimensional approximation of this set. In order to maintain technicalities to a minimum, we consider only quadratic interactions. The Dirichlet energy of such a set is defined as
$$
F(u)=\sum_{(i,j)\in\mathcal E}(u(i)-u(j))^2
$$
on functions $u:X\to \mathbb R$. A discrete-to-continuum and dimensionally reduced limit is then obtained by considering a scaled version of the energies
$$
F_\e(u)=\sum_{(i,j)\in\mathcal E}\e^{d-2} (u_i-u_j)^2
$$
defined on functions $u:\e X\to \mathbb R$, where we use the notation $u_i=u(\e i)$, and taking their limit in a suitable sense as $\e\to0$. Note that we may interpret
$$
\e^{-1}(u_i-u_j) = |i-j|\Bigl({u_i-u_j\over \e|i-j|}\Bigr)
$$
as an inhomogeneous difference quotient, so that $F_\e$ represent discrete versions of an (inhomogeneous) Dirichlet integral, whose general continuous counterpart is of the form
\begin{equation}\label{conttf}
{1\over \e^k}\int_{\e E}f(\nabla u)\,dx,
\end{equation} with $E$ a subset of $\R^d\times \R^k$ uniformly bounded in the last $k$ variables. Energies of the form \eqref{conttf} are the prototype of thin-structure energies on the continuum (see e.g.~\cite{Brai02,Handbook}), which have been treated extensively in the last thirty years. Among the many contribution to the subject we recall the seminal paper by Le Dret and Raoult \cite{LDR} which gives a general dimension-reduction formula  when $E=\R^d\times [0,1]$ through a lower-dimensional quasiconvexification process.  Moreover, a general compactness and integral-representation theorem has been proved by Braides, Fonseca and Francfort \cite{BFF}, which interpret lower-dimensional quasiconvexification through a homogenization formula, and extend the analysis to general thin films with varying profiles.  In their approach they deal with $E$ that can be seen as a subgraph of a function defined on $\R^d$. In our case, even if a continuum set $E$  corresponding to $X$ can be constructed, it may not be a subgraph, as it might have holes or even possess a more complex topology. Note that the assumption that the integration be performed on the scaled $\e E$ cannot be extended to arbitrary $E_\e$ since in that case the limit might not be simply $d$-dimensional if the complexity of the topology increased as $\e\to0$ (see the example by Braides and Bhattacharya \cite{BBa}). 
Finally, we note that asymptotic analysis of thin objects can be interpreted as an intermediate step in the study of structures with very fast oscillating profile \cite{BCZ} (see also \cite{AB}, and e.g.~\cite{Gau} for an example of application in a continuum geometry). 

\smallskip
Discrete-to-continuum analyses for lattice energies are usually performed after identification of functions defined on (portions of) lattices with their piecewise-constant interpolations. This identification allows to embed families of energies in a common Lebesgue-space environment (see the seminal paper by Alicandro and Cicalese \cite{AC04}). Using this approach, a discrete-to-continuum analog for thin films of the Braides, Fonseca and Francfort theory, has been studied by Alicandro, Braides and Cicalese \cite{ABC} (see also \cite{Sc}, and the work \cite{BRS} for a connection with aperiodic lattices). Due to the great generality of our discrete set $X$, we will not directly extend functions defined on $\e X$ but follow a dimension-lowering coarse-graining approach: to each function $u_\e:\e X\to \R$ we associate the function $\overline u_\e:\e T\Z^d\to \R$ where $\overline u_\e(\e Tl)$ is obtained by averaging $u_\e$ on $X\cap \bigl((\e Tl+ \{0,\ldots, T-1\}^d)\times \{0,\ldots, M-1\}^k\bigr)$. We then prove that energy bounds on $u_\e$ imply that the piecewise-constant interpolations of the corresponding $\overline u_\e$ are precompact in $L^2_{\rm loc}(\R^d)$ and their limit is in $H^1_{\rm loc}(\R^d)$. In this way a dimensionally reduced continuum parameter can be defined. In order to relate the original $u_\e$ to this limit, a Poincar\'e inequality must be used at scale $\e$, which shows that the original $u_\e$ converge to $u$ in a `perforated domain' fashion (see e.g.~\cite{BCP-doublepor}). Both the coarse-graining and the Poincar\'e-type inequality are very reminiscent of the $p$-connectedness approach by Zhikov \cite{Z}, and of its use in the homogenization of singular structures by Braides and Chiad\`o Piat \cite{BCP-ss}, even though in those papers $p$-connectedness is stated fo local functionals depending on the gradient. Here we deal with non-local interactions, even though the non-locality weakens as $\e\to0$, and some additional care has to be taken, similarly to the case of the homogenization of convolution-type energies (see \cite{AABPT, BCD, BPiat}).

\smallskip
The paper is organized as follows. In Section \ref{Sett} we introduce the notation for the environment $X\subset \R^d\times \R^k$ and for the energies that we consider, which are a little more general than those described above in that a more general inhomogeneity is allowed (introducing interactions coefficients $a_{ij}$) and the energies are localized by considering interactions parameterized on a set $\Omega\subset \R^d$. Section \ref{Twoconn} is devoted to the definition of coarse-grained functions, and to the statement and proof of the two-connectedness property and of a Poincar\'{e}-Wirtinger's inequality. Section \ref{Comp} contains a compactness result for coarse-grained functions whose proof relies on the two-connectedness property and the corresponding convergence of the original functions. The limit defines a function on a subset of $\R^d$. Section~\ref{Bounda} contains a result that allows to consider boundary-values on `lateral boundaries' of thin films. A homogenization theorem for quadratic energies defined on $\e X$ is stated in Section~\ref{Homsec}. Its proof is subdivided into a lower bound by blow-up and an upper bound by a direct construction. Moreover, an application to the description of the asymptotic behaviour of boundary-value problems is also described. Finally, Section \ref{Examples} contains some simple examples illustrating some possible non-trivial shapes of the thin structures we consider.

\paragraph*{Notation}
\begin{itemize}
\item The letter $C$ denotes a generic strictly positive constant not depending on the parameters of the problem considered, whose value may be different at every its appearance. 

\item If $x, y\in \mathbb R^d$ then $x\cdot y$ denotes their scalar product. If $t\in \mathbb R$ then $\lfloor t\rfloor$ is its integer part.

\item For $T\in\mathbb{N}$, we denote by $\cubeTd$ the $d$-dimensional semi-open cube of side length $T$; {\it i.e.}, $\cubeTd:= [0,T)^d$. If $T=1$, we simply write $Q_d= Q_{1,d}$. For $l\in\Z^d$, $\lcubeTd:= lT + [0,T)^d$ and for $T=1$, we write $Q^l_d=Q^l_{1,d}$.

\item $\cubeTk$ denotes the $k$-dimensional semi-open cube of side length $T$; {\it i.e.}, $\cubeTk:= [0, T)^k$. For $T=1$, we set $Q_k=Q_{1,k}$ and $Q^n_k=n+ Q_k$ if $n\in\Z^k$.

\item For any measurable set $\Omega$ and $u\in L^1(\Omega)$,  $\fint_{\Omega} u(x)dx$ denotes the average of $u$ on $\Omega$; {\it i.e.},
   \begin{equation}
   \notag
   \fint_{\Omega}u(x)dx := {1\over |\Omega|}\int_{\Omega}u(x)dx,
   \end{equation}
where $|\cdot|$ stands for the Lebesgue measure.

\item For any open set $\Omega\in\R^d$ and for any $\delta>0$,  we let $
    \Omega(\delta):=\{x^d\in\Omega\hspace{0.03cm}:\hspace{0.03cm} \dist(x, \partial \Omega)>\delta   \}
$.
\end{itemize}

\section{Setting of the problem}\label{Sett}
In the following  $X$ will be a fixed subset of $\Z^d\times\{0,\dots, T-1\}^k$, with $d,k\geq 1$ and $T\in\mathbb{N}$. We assume that 
    \begin{itemize}
    \item[(i) ] $X$ is $T$-periodic in $e_1,\dots, e_d$;
    \item[(ii) ]  $X$ is connected in the following sense: there exists $\mathcal{E}\subset X\times X$ such that for all $i,j\in X$ there exists a sequence $\{i_n\}_{n=0}^{N}$ of points of $X$, with $i_0=i$ and $i_N=j$, such that the segment $(i_n, i_{n+1})\in\mathcal{E}$. Moreover, the set $\mathcal{E}$ is $T$-periodic; {\it i.e.}, if the segment $(i, j)$ belongs to $\mathcal{E}$, then, for any $m=1,\dots, d$, the segment $(i+ Te_m, j+Te_m)$ belongs to $\mathcal{E}$;
    \item[(iii)] the set $\mathcal{E}$ is equi-bounded; {\it i.e.}, there exists $R>0$ such that 
                    \begin{equation*}
                    \max\{|i-j|\hspace{0.03cm}:\hspace{0.03cm} (i,j)\in \mathcal{E} \}\leq R.
                    \end{equation*}
     \end{itemize} 

   {\rm Note that it is not restrictive to assume that $R\leq T$, upon taking a larger period.}
\smallskip

Let $a_{ij}$ be $T$-periodic coefficients in $e_1,\dots, e_d$; {\em i.e.}, 
$$
a_{i+Te_m\,j+T e_m}= a_{ij} \hbox{ for all } i,j \in  X,\  m\in\{1,\ldots, d\},
$$
such that $a_{ij}>0$ if $(i,j)\in\mathcal{E}$ and $a_{ij}=0$ if $(i,j)\notin\mathcal{E}$. For $\e>0$, we introduce the family of functionals $F_\e$ defined on functions $u:\e X\to \R$ by
          \begin{align}
          \notag
          F_\e(u) := \sum_{i,j\in X} \e^d a_{ij} \left({u_i-u_j\over \e} \right)^2, %\label{energy}
          \end{align}
where $u_i:=u(\e i)$. Note that also the case $a_{ij}=1$  if $(i,j)\in\mathcal{E}$ is non trivial. Note moreover that, by the periodicity of $a_{ij}$,  there exists a positive constant $C$ such that $C\le a_{ij}\le 1/C$ if $(i,j)\in\mathcal{E}$, so that $F_\e$ is estimated from above and below by the energy corresponding to $a_{ij}=1$  if $(i,j)\in\mathcal{E}$.

More in general, we will consider `localized' versions of energies $F_\e$, limiting interactions to $i,j\in X$ such that $\e i, \e j \in \Omega \times \e\cubeTk$ for some Lipschitz open subset $\Omega$ of $\R^d$.

   \begin{remark}
  {\rm In the notation above, we can include also the case of 
       \begin{equation*}
       X\subset \Z^d\times\prod_{n=1}^{k}\{0, \dots, M_n-1\},
       \end{equation*}
with $T_m\geq 1$, $m=1,\dots, d$, and $M_n\geq 1$, $n=1,\dots, k$,  and $X$ $T_m$-periodic in $e_m$, for any $m=1,\dots, d$. In this case,  we take $T={\rm l.c.m.}\{ T_1,\dots, T_d, M_1,\dots, M_k\}$.}
   \end{remark}

%   \begin{example}
%   \begin{itemize}
 %    \begin{figure}[h]
  %     \centering
   %    \includegraphics[scale=0.30]{DNA_Gray.png}
    %   \end{figure}
   %\end{itemize}
%   \end{example}

\section{Two-connectedness and Poincar\'{e}-Wirtinger's inequality}\label{Twoconn}
In this section we prove two technical lemmas, which will allow to use some compactness results for systems of nearest-neighbour interactions. To that end, we define a coarse-grained lower-dimensional variable as follows.
Let $u$ be a real-valued function defined on $X$. For $l\in\Z^d$, set
   \begin{equation}
   \label{meanvalue}
   \tilde{u}^l := {1\over \#[(\cubeTd\times \cubeTk)\cap X] }\sum_{i\in (\lcubeTd\times \cubeTk)\cap X} u_i, 
   \end{equation}
  with $u_i= u(i)$.

  The first result of this section states that a nearest-neighbour interaction energy on the coarse-grained variable is (locally) dominated by the energy on $X$.
  \begin{prop}
  There exist $C=C(X)>0$ and $M>0$ such that  
    \begin{align}
    |\tilde{u}^l-\tilde{u}^{l'}|^2\leq C \sum_{i,j\in [(\cubeTd^l\cup \cubeTd^{l'}+(-M, M)^d)\times\cubeTk]\cap X } |u_i-u_j|^2, \label{2connectedness}
    \end{align}
    for any $l, l'\in \Z^d$ such that $|l-l'|=1$.
  \end{prop}
  \begin{proof}
  Using definition \eqref{meanvalue} of $\tilde{u}^l$ and the change of indices $j=i+Te_m$, for some $m=1,\dots, d$, combined with the H\"{o}lder inequality, we deduce that 
        \begin{align}
        |\tilde{u}^l-\tilde{u}^{l'}|^2&= {1\over[\#((\cubeTd\times \cubeTk)\cap X )]^2}\biggl|\sum_{i\in (\cubeTd^l\times \cubeTk)\cap X }u_i-\sum_{j\in(\cubeTd^{l'}\times \cubeTk)\cap X }u_j  \biggr|^2\notag\\
        &={1\over[\#((\cubeTd\times \cubeTk)\cap X )]^2}\biggl|\sum_{i\in (\cubeTd^{l}\times \cubeTk)\cap X }(u_i-u_{i+Te_m})  \biggr|^2\notag\\
        &\leq {1\over\#[(\cubeTd\times \cubeTk)\cap X ]}\sum_{i\in(\cubeTd^l\times \cubeTk)\cap X }|u_i-u_{i+Te_m}|^2\label{est1}.
        \end{align}
  The connectedness of $X$ ensures that for all $i\in(\cubeTd^l\times\cubeTk)\cap X$, there exists a sequence $\{j_n\}_{n=0}^{N_i}$ of points in $X$ with $j_0=i$ and $j_{N_i}=i+Te_m$ such that $(j_n, j_{n+1})\in\mathcal{E}$. Let $\gamma$ be path joining $i$ and $i+Te_m$ through the points $j_1,\dots, j_{N_i-1}$. Such a path $\gamma$ is contained in $[(\cubeTd^l\cup\cubeTd^{l'} +(-M, M)^d)\times\cubeTk]\cap X$, for some $M>0$ large enough independent of the point $i$. For any $i\in(\cubeTd^l\times\cubeTk)\cap X$, we write that 
        \begin{align}
        u_i-u_{i+Te_m}=\sum_{n=1}^{N_i} (u_{j_{n-1}}-u_{j_n}) \notag,
        \end{align}
 so that,  due to \eqref{est1} combined with the H\"{o}lder inequality, we have that 
        \begin{align}
         |\tilde{u}^l-\tilde{u}^{l'}|^2 &\leq  {1\over\#[(\cubeTd\times \cubeTk)\cap X ]}\sum_{i\in (\cubeTd^l\times \cubeTk)\cap X }\biggl|\sum_{n=1}^{N_i} (u_{j_{n-1}}-u_{j_n})\biggr|^2\notag\\
         &\leq {1\over\#[(\cubeTd\times \cubeTk)\cap X ]} \sum_{i\in (\cubeTd^l\times \cubeTk)\cap X } N_i\sum_{n=1}^{N_i}|u_{j_{n-1}}-u_{j_n}|^2\notag\\
         &\leq {\max\{N_i\hspace{0.03cm}:\hspace{0.03cm} i\in (\cubeTd^l\times\cubeTk)\cap X \}\over\#[(\cubeTd\times\cubeTk)\cap X ]}\sum_{i,j\in [(\cubeTd^l\cup\cubeTd^{l'}+(-M, M)^d)\times \cubeTk]\cap X } |u_i-u_j|^2,\notag
        \end{align}
 where in the last inequality we have used the fact that   $[(\cubeTd^l\cup\cubeTd^{l'}+(-M, M)^d)\times \cubeTk]\cap X$ contains the  path $\gamma$ joining $i$ and $i+Te_m$ for all $i\in (\cubeTd^l\times\cubeTk)\cap X$. %Finally, using the property of $a_{ij}$ given by \eqref{aij}, we conclude that 
      % \begin{align}
      % |\tilde{u}^l-\tilde{u}^{l'}|^2 &
      % &\leq {1\over C}{\max\{N_i\hspace{0.03cm}:\hspace{0.03cm} i\in \cubeTd^l\times\cubeTk\cap X \}\over\#(\cubeTd\times\cubeTk\cap X )}\sum_{i,j\in (\cubeTd^l\cup\cubeTd^{l'}+(-M, M)^d)\times \cubeTk\cap X } a_{ij} |u(i)-u(j)|^2,\notag
      % \end{align}
This proves the desired inequality.
  \end{proof}
  \begin{remark}\rm
  In order to reduce the number of parameters, we can choose $M=T$, up to substituting $T$ with a multiple and taking a slightly larger $M$.
  \end{remark}
We point out that in the following \eqref{meanvalue} and \eqref{2connectedness} will be applied to functions $u:\e X\to\R$, where $u_i$ stands for $u(\e i)$ as in the notation introduced above.

\begin{comment} %2-CONNECTEDNESS WITH $u_\e$
For $l\in\Z^d$, we define $\tilde{u}^l_\e$ as
    \begin{equation}
    \label{deftildeule}
    \tilde{u}^l_\e:= {1\over \#(\cubeTd\times \cubeTk\cap X)}\sum_{i\in\cubeTd^l\times \cubeTk\cap X}u_i,
    \end{equation}
A consequence of \eqref{2connectedness} is the rescaled $2$-connectedness inequality  given by
 \begin{align}
    |\tilde{u}^l_\e-\tilde{u}^{l'}_\e|^2\leq C \sum_{i,j\in (Q^l_T\cup Q^{l'}_T+(-M, M)^d)\times\cubeTk\cap X}|u_i-u_j|^2, \label{rescaled2connect}
    \end{align}
for any $l, l'\in \Z^d$ such that $|l-l'|=1$. %Indeed, from \eqref{2connectedness}, it follows that
%    \begin{align}
 %    |\tilde{u}^l_\e-\tilde{u}^{l'}_\e|^2&=\biggl|{1\over \#(\cubeTd\times\cubeTk\cap X) }\biggl( \sum_{i\in\cubeTd^l\times \cubeTk\cap X}u_i-\sum_{i\in\cubeTd^{l'}\times \cubeTk\cap X}u_i \biggr) \biggr|^2\notag\\
  %   &\leq \e^2C \sum_{i,j\in((\cubeTd^l\cup\cubeTd^{l'}+(-M, M)^d))\times \cubeTk\cap X} |u_i - u_j|^2\notag
     %&=\e^{2-d}C \sum_{i,j\in(\e(\cubeTd^l\cup\cubeTd^{l'}+(-M, M)^d))\times \cubeTk\cap X_\e} |u( i^d, \e i^k) - u( j^d, \e j^k)|^2,\notag
 %   \end{align}
 %which proves \eqref{rescaled2connect}.
\end{comment}  

\smallskip
 Now, we show a Poincar\'{e}-Wirtinger inequality. This will be used to recover information on the original functions $u$ from their coarse-grained versions.
      \begin{prop}
      \begin{itemize}
       \item[\rm(i)] There exists $C=C(X)>0$ such that, for any $l\in\Z^d$,  
         \begin{equation}
          \sum_{i\in (\cubeTd^l\times\cubeTk)\cap X}|u_i-\tilde{u}^l|^2\leq C \sum_{i,j\in (\cubeTd^l\times\cubeTk)\cap X}|u_i-u_j|^2;\notag
          \end{equation}
      \item[\rm(ii)] there exist positive constants $C$ and $M$ such that, for any $l\in\Z^d$, 
           \begin{equation}
           \sum_{i\in (\cubeTd^l\times\cubeTk)\cap X}|u_i-\tilde{u}^l|^2\leq C \sum_{i,j\in [(\cubeTd^l+(-M,M)^d)\times\cubeTk]\cap X}a_{ij}|u_i-u_j|^2.\label{rescaledPW}
           \end{equation}
      \end{itemize}
      \end{prop}
      \begin{proof}
(i)
      Using  definition \eqref{meanvalue} of $\tilde{u}^l$ and thanks to the H\"{o}lder inequality, we deduce that 
          \begin{align}
                   \sum_{i\in (\cubeTd^l\times\cubeTk)\cap X}|u_i-\tilde{u}^l|^2&= {1\over [\#((\cubeTd\times\cubeTk)\cap X)]^2}\sum_{i\in (\cubeTd^l\times\cubeTk)\cap X}\biggl| \sum_{j\in (\cubeTd^l\times\cubeTk)\cap X}(u_i-u_j) \biggr|^2\notag\\
                   &\leq {1\over \#[(\cubeTd\times\cubeTk)\cap X]}\sum_{i,j\in (\cubeTd^l\times\cubeTk)\cap X}|u_i-u_j|^2\label{estpw-withoutaij} 
                   %&\leq {1\over C} {1\over \#(\cubeTd\times\cubeTk\cap X)}\sum_{i,j\in \cubeTd^l\times\cubeTk\cap X}a_{ij}|u(i)-u(j)|^2\notag,
          \end{align}
      which concludes the proof.
      
(ii)
Since $X$ is connected and due to the boundedness and periodicity properties of the coefficient $a_{ij}$, there exists $M>0$ large enough such that if $i,j\in (\lcubeTd\times \cubeTk)\cap X$, then there exists a path $\gamma$ joining $i$ and $j$ which is contained in $[(\lcubeTd+(-M, M)^d)\times\cubeTk]\cap X$. From \eqref{estpw-withoutaij}, we deduce \eqref{rescaledPW} as desired.
      \end{proof}

\begin{comment}
Repeating similar arguments as \eqref{rescaled2connect}, the rescaled Poincar\'{e}-Wirtinger inequality is given by
    \begin{align}
     \sum_{i\in \lcubeTd\times\cubeTk\cap X}|u_i-\tilde{u}^l_\e|^2\leq C \sum_{i,j\in \lcubeTd\times\cubeTk\cap X}|u_i-u_j|^2.\notag
    \end{align}
\end{comment}

\section{A compactness result}\label{Comp}
In this section, we are going to show that sequences with equi-bounded energy are compact in $L^2$ with limit in $H^1_{\loc}(\R^d)$.\par Let $\Omega$ be an open set of $\R^d$ with Lipschitz boundary. For $\e>0$, let $u_\e$ be a family of functions $u_\e: (\Omega\times \e\cubeTk)\cap\e X\to\R$. Setting $ I_\e=I_\e(\Omega):= \{l\in\Z^d\hspace{0.03cm}:\hspace{0.03cm}\e\lcubeTd\subset\Omega\}$, we define a  piecewise-constant function $\overline{u}_\e$ 
in $L^2(\Omega)$  by
     \begin{align}
     \overline{u}_\e(x^d):=\sum_{l\in I_\e}\tilde{u}^l_\e\chi_{\e \lcubeTd}(x^d),\label{defoverlineue}
     \end{align}
where $\tilde{u}^l_\e$ is given by \eqref{meanvalue}, with $u_i=u^\e_i:= u_\e(\e i)$, and $\chi_{\e\lcubeTd}$ is the characteristic function of the cube $\e\lcubeTd$.
\par In order to introduce the convergence of $u_\e$,  we need to identify  real-valued functions $u$ defined on $\e X$ with piecewise-constant interpolations as follows. We introduce the set
    \begin{align}
    \mathcal{C}_\e(\Omega):= \Bigl\{& u:\R^d\times \cubeTk\to\R\hspace{0.03cm}:\hspace{0.03cm} u \mbox{ is constant on } \e Q^l_d \times \e Q^{n}_k\notag\\
   & \hspace*{1cm} \mbox{ for } (l,n)\in \Bigl(\Z^d\cap{1\over\e}\Omega\Bigr)\times\{0,\dots,T-1\}^k \Bigr\}, \label{defCeOmegaX}
    \end{align}
so that a  function $u:(\e \Z^d\cap\Omega)\times\e\{0,\dots,T-1\}^k\cap \e X\to\R$ can be identified with its extension belonging to $\mathcal{C}_\e(\Omega)$. We say that the family of function $u_\e$ in $\mathcal{C}_\e(\Omega)$ converges to $u\in H^1(\Omega)$  if
         \begin{equation}
          \label{defconvergence}
         \overline{u}_\e \to u \quad\mbox{in }L^2_\loc(\Omega).
         \end{equation}
From this convergence, we will obtain that 
     \begin{equation}\label{convsuX}
     \int_{\widetilde{\Omega}_\e}|u_\e -u|^2\chi_{\displaystyle(\cup_{i\in X}\e Q^{i^d}_d\times Q_k^{i^k}  )}\to 0,
     \end{equation}
where we have set $i=(i^d, i^k)\in\Z^d\times\{0,\dots, T-1\}^k$ and  $\widetilde{\Omega}_\e$ is given by
    \begin{equation}
    \label{deftildeomegae}
    \widetilde{\Omega}_\e := \bigcup_{l\in{\cal L}_\e}\e\lcubeTd\times\cubeTk,
    \end{equation}
and ${\cal L}_\e:= \{l\in Z^d: \dist(\e l, \partial \Omega)>2\e\sqrt{d}T   \}$. The next proposition provides a compactness result for $\overline{u}_\e$ using the analysis of nearest-neighbour interactions in  \cite{AC04}.%whose proof follows the arguments of \cite[Lemma 2.1]{BCP08}.
 \begin{prop}[Compactness]
 \label{prop:cptresultmeanvalue}
Let $\Omega$ be an open set of $\R^d$ with Lipschitz boundary. Let $u_\e$ be a family of functions defined on $(\Omega\times \e\cubeTk)\cap\e X$ such that 
     \begin{equation}
     \sum_{l\in\Z^d}\sum_{i\in(\lcubeTd\times\cubeTk)\cap X}\e^d|u^\e_i|^2\leq C, \qquad \forall\e>0,\label{equibounormue}
     \end{equation}
 where $u^\e_i=0$ if $i\notin[({1\over\e}\Omega\cap\lcubeTd)\times\cubeTk]\cap X$.
 and 
      \begin{equation}
      \sum_{l,l'\in {\cal L}_\e, |l-l'|=1}\sum_{i,j\in[(\cubeTd^l\cup\cubeTd^{l'}+(-M, M)^d) \times\cubeTk]\cap X}\e^{d-2} |u^\e_i - u^\e_j|^2\leq C, \qquad\forall\e>0. \label{equiboungrad}
      \end{equation}
  Then, up to a subsequence, the family  $\overline{u}_\e$, given by \eqref{defoverlineue}, strongly converges in  $L^2_{\loc}(\Omega)$ to some $u\in H^1(\Omega)$.
 \end{prop}
 \begin{proof} First, we show that $\overline{u}_\e$ weakly converges in $L^2_\loc$ to some $u$. Indeed, from \eqref{equibounormue}, we deduce that the norm $\|\overline{u}_\e\|_{L^2(\Omega)}$ is bounded which implies the weak convergence of $\overline{u}_\e$. Morever, thanks to assumption \eqref{equiboungrad}, an application of  \cite[Proposition 3.4]{AC04} provides us that $u\in H^1(\Omega)$.\par Now, we prove the strong convergence in $L^2_\loc(\Omega)$ of $\overline{u}_\e$. To this end, 
 the key tool is the Compactness Criterion by Fr\'{e}chet and Kolmogorov (see,  e.g. \cite[Theorem 4.26]{Bre10}). In other words, we have to prove that, for any $\Omega''\subset\subset\Omega'$ and for any $\eta>0$, there exists $\delta>0$, with $\dist(\Omega'', \R^d\setminus\Omega')>\delta$, such that for every $h\in\R^d$, with $|h|<\delta$, then
       \begin{equation}
       \|\tau_h\overline{u}_\e-\overline{u}_\e\|_{L^2(\Omega'')}<\eta,\label{estcptthmFK}
       \end{equation}
  where $\tau_h\overline{u}_\e(x):= \overline{u}_\e(x+h)$.  Assume that $h=\lambda e_m$, for some $m=1\dots, d$. The inequality \eqref{estcptthmFK} for every $h\in\R^d$ is obtained by triangle inequality. Fix $\Omega''\subset\subset\Omega'$ and set
      \begin{equation}
      \mathcal{I}_\e:=\{l\in{\cal L}_\e\hspace{0.03cm}:\hspace{0.03cm} \Omega'\subset \cup_{l}\e\lcubeTd\subset\Omega \}.\notag
      \end{equation}
   Take $x\in\e\lcubeTd$. Hence, we have that $x\in\e\lcubeTd$ and $(x+h)\in\e\cubeTd^{l'}$, for some $l,l'\in\mathcal{I}_\e$. By definition of $\overline{u}_\e$ given by \eqref{defoverlineue}, we deduce that 
       \begin{align}
       |\tau_h\overline{u}_\e(x)-\overline{u}_\e(x)|^2=|\overline{u}_\e(x+h)-\overline{u}_\e(x)|^2= |\tilde{u}^l_\e-\tilde{u}^{l'}_\e|\label{tauee-ue}
       \end{align}
    Since $l$ and $l'$ are not necessarily such that $|l-l'|=1$, we need to re-write the two-connectedness inequality in terms of  non-neighbouring cubes. In order to show this, let  $S_{ll'}$ be union of neighbouring cubes joining $\e\lcubeTd$ and $\e\cubeTd^{l'}$ such that each two consecutive cubes have one face in common; {\it i.e.}, $S_{ll'}=\bigcup_{n=0}^{N_\e}\e\cubeTd^{l_n}$, with $|l_n-l_{n-1}|=1$, $l_0=l$ and $l_{N_\e}=l'$. Note that the number $N_\e$ of the cubes $\e\lcubeTd$ contained in stripes of cubes joining $\e\cubeTd^l$ and  $\e\cubeTd^{l'}$ is of order $|h|T^{-1}\e^{-1}$. Hence, thanks to inequality \eqref{2connectedness}, we deduce that
        \begin{align}
        |\tilde{u}^l_\e-\tilde{u}^{l'}_\e|^2&= \biggl|\sum_{n=1}^{N_\e}(\tilde{u}^{l_n}_\e-\tilde{u}^{l_{n-1}}_\e ) \biggr|^2\notag\\
        &\leq |h|T^{-1}\e^{-1}\sum_{n=1}^{N_\e}|\tilde{u}^{l_n}_\e-\tilde{u}^{l_{n-1}}_\e|^2\notag\\
        &\leq C |h|T^{-1}\e^{-1} \sum_{n=1}^{N} \sum_{i,j\in[(\cubeTd^{l_n}\cup\cubeTd^{l_{n-1}}+(-M, M)^d)\times\cubeTk]\cap X} |u^\e_i - u^\e_j|^2\notag\\
        &\leq C|h|T^{-1} \e^{-1} \sum_{i,j\in [(S_{ll'}+(-M, M)^d )\times \cubeTk]\cap X}|u^\e_i - u^\e_j|^2.\notag
        \end{align}
  Plugging the above inequality in \eqref{tauee-ue}, we obtain that 
       \begin{align*}
       |\tau_h\overline{u}_\e(x)- \overline{u}_\e(x)|^2 &\leq C|h|T^{-1}\e^{-1}\sum_{i,j\in [(S_{ll'}+(-M, M)^d )\times \cubeTk]\cap X}|u^\e_i - u^\e_j|^2\notag\\
       &= C|h|T^{-1}\e^{-1}\sum_{i,j\in (S(x,h)\times \cubeTk)\cap X}|u^\e_i - u^\e_j|^2,\notag
       \end{align*}
   where we have set $S(x,h):= S_{ll'}+(-M, M)^d$ which depends on $x$ and $h$. Now, an integration  with respect to $x\in\e\lcubeTd$ yields
       \begin{align}
       \int_{\e\lcubeTd}|\tau_h\overline{u}_\e(x)- \overline{u}_\e(x)|^2 dx&\leq C|h|T^{-1}\e^{-1}  \int_{\e\lcubeTd}\biggl(\sum_{i,j\in (S(x,h)\times \cubeTk)\cap X}|u^\e_i - u^\e_j|^2\biggr)dx\notag\\
       &\leq C|h|T^{d-1}\e^{d-1} \sum_{i,j\in (S(\e\lcubeTd, h)\times\cubeTk)\cap X}|u^\e_i - u^\e_j|^2,\notag
       \end{align}
    where we have used the fact that 
        \begin{equation}
        S(x,h)\subset S(\e\lcubeTd,h):=\bigcup\{S(x,h)\hspace{0.03cm}: \hspace{0.03cm} x\in\e\lcubeTd\}.\notag
        \end{equation}
  Summing over $\mathcal{I}_\e$, we have that 
      \begin{align*}
      \sum_{l\in\mathcal{I}_\e}\int_{\e\lcubeTd}|\tau_h\overline{u}_\e(x)- \overline{u}_\e(x)|^2 dx &\leq C|h|T^{d-1}\e^{d-1}  \sum_{l\in\mathcal{I}_\e}\sum_{i,j\in (S(\e\lcubeTd, h)\times\cubeTk)\cap X}|u^\e_i - u^\e_j|^2,\notag
      \end{align*}
   which implies that
       \begin{align}
       \int_{\Omega''}|\tau_h\overline{u}_\e(x)- \overline{u}_\e(x)|^2 dx \leq \sum_{l\in\mathcal{I}_\e}\int_{\e\lcubeTd}|\tau_h\overline{u}_\e(x)- \overline{u}_\e(x)|^2 dx \leq C |h|, \notag
       \end{align}
       where we have used assumption \eqref{equiboungrad} and the fact that the number of indices $l$ and $l'$ such that $S(\e\lcubeTd, h)\cap S(\e\cubeTd^{l'},h)\neq\emptyset $ is of the order $|h|T^{-1}\e^{-1}$ (the ratio between the size of $S(\e\cubeTd^l,h)$ and the size of $\e\cubeTd^l$). This conclude the proof of \eqref{estcptthmFK}. \par Finally, applying the compactness criterion, it follows that, up to a subsequence, $\overline{u}_\e\to v.$ Since we already know that $\overline{u}_\e\rightharpoonup u$, we conclude that $v=u$, which is the desired claim.\end{proof}

The next proposition provides a convergence result in the sense of \eqref{convsuX}.

 \begin{prop}
 Let $\Omega$ be an open set of $\R^d$ with Lipschitz boundary. Let $u_\e$ be a sequence of functions defined on $\e X$ such that
      \begin{equation}
      \label{equiboundcpt}
      \sup_{\e>0}\biggl( \sum_{l\in{\cal L}_\e}\sum_{i\in(\lcubeTd\times\cubeTk)\cap X}\e^d|u^\e_i|^2+ F_\e(u_\e) \biggr)\leq C.    
      \end{equation}
   Then, up to a subsequence, we have that
     \begin{equation}
     \notag
      \int_{\widetilde{\Omega}_\e}|u_\e-u|^2\chi_{\cup_{i\in X}\e Q_d^{i^d}\times\cubeTk^{i^k}}\to 0 
      \end{equation}
  where $u\in H^1(\Omega)$ is the strong limit in $L^2_{\loc}(\Omega)$ of the sequence $\overline{u}_\e$ and $\widetilde{\Omega}_\e$ is given by \eqref{deftildeomegae}.
 \end{prop}
 \begin{proof}
  Set $x=(x^d, x^k)\in\e\lcubeTd\times\cubeTk$ and recall that $u_\e$ is defined on $\e\lcubeTd\times\e\cubeTk$. Hence, 
     \begin{align}
     \int_{\tilde{\Omega}_\e} |(u_\e-u)(x)\chi_{\displaystyle\cup_{i\in X}\e Q_d^{i^d}\times Q_k^{i^k}}(x)|^2dx&\leq\int_{\tilde{\Omega}_\e} |(u_\e-\overline{u}_\e)(x)\chi_{\displaystyle\cup_{i\in X}\e Q_d^{i^d}\times Q_k^{i^k}}(x)|^2dx\notag\\
     &\quad +\int_{\tilde{\Omega}_\e} |(\overline{u}_\e-u)(x)\chi_{\displaystyle\cup_{i\in X}\e Q_d^{i^d}\times Q_k^{i^k}}(x)|^2dx.\label{f0}
     %&\leq\sum_{i\in X}\int_{\e Q^{i^d}_d\times Q^{i^k}_k}|u_\e(x^d, \e x^k) - \overline{u}_\e(x^d)|^2dx^ddx^k\notag\\ &\quad+\sum_{i\in X}\int_{\e Q^{i^d}_d\times Q^{i^k}_k}|\overline{u}_\e(x^d)-u(x^d)|^2dx^ddx^k \notag\\
     %&\leq  \sum_{i\in X}\int_{\e Q^{i^d}_d\times Q^{i^k}_k}|u_\e(x^d, \e x^k) - \overline{u}_\e(x^d)|^2dx^ddx^k\notag\\
     %&\quad + T^d\int_{\color{red}\Omega}|\overline{u}_\e(x^d) - u(x^d)|^2dx^d.
     \end{align}
    From Proposition \ref{prop:cptresultmeanvalue}, we know that $\overline{u}_\e$ strongly converges to $u$ in $L^2_\loc(\Omega)$, so that the second integral in \eqref{f0} vanishes as $\e\to 0$. \par In order to estimate the first integral of \eqref{f0}, the key tool is the  Poincar\'{e}-Wirtinger inequality given by \eqref{rescaledPW}. Indeed, due to the fact that $u_\e$ is constant on $\e Q^l_d\times Q_k$  and $\overline{u}_\e$ is constant on $\e\lcubeTd\times Q_k$, we deduce that
         \begin{align}
         \int_{\tilde{\Omega}_\e} |(u_\e-\overline{u}_\e)(x)\chi_{\displaystyle\cup_{i\in X}\e Q_d^{i^d}\times Q_k^{i^k}}(x)|^2dx& = \e^d\sum_{l\in{\cal L}_\e}\sum_{i\in (\lcubeTd\times\cubeTk)\cap X } |u^\e_i-\tilde{u}^l_\e|^2\notag\\
         &\leq \e^d \sum_{l\in{\cal L}_\e} \sum_{i,j\in[(\cubeTd+(-T, T)^d)\times\cubeTk]\cap X} a_{ij} |u^\e_i-u^\e_j|^2\notag\\
         &\leq \e^2 F_\e(u_\e).\notag
         \end{align}
    From this, combined with assumption \eqref{equiboundcpt}, we have that also the first integral of \eqref{f0} goes to $0$ as $\e\to 0$, which concludes the proof.
 \end{proof}

\section{Treatment of boundary data}\label{Bounda}
In this section we prove a classical lemma which allows to match boundary conditions. For future reference we prove it  in a general form. \par For any $u\in H^1(\Omega)$, we define the sequence $v_\e$ on $\e X$  by
    \begin{equation}
    \label{def:ve}
    v^\e_i = v_\e(\e i):= \fint _{\e i^d+\e Q_{d}} u(x)dx.
    \end{equation}
We have that $v_\e$ converges to $u$ with respect to convergence \eqref{defconvergence}.
For any bounded open set $A$ and for $\delta>0$, we define  
$A(\delta):=\{x\in A\hspace{0.03cm}:\hspace{0.03cm} \dist(x, \partial A)>\delta \}$.
  \begin{lemma}
  \label{lemma:treatmentBC}
  Let $A$ be a bounded and open set of $\Omega$ with Lipschitz boundary. Let $u_\e$ be a sequence converging to $u\in H^1(\Omega)$ with respect to convergence \eqref{defconvergence}. For any $\delta>0$, there exists a sequence $w_\e$ converging to $u$ with respect convergence \eqref{defconvergence} such that 
      \begin{align}
      \notag
      w_\e &= u_\e, \qquad\mbox{if } i\in (A(2\delta)\times\cubeTk)\cap X,\\
      \notag
      w_\e&=v_\e, \qquad\mbox{if } i\in (A\setminus A(\delta)\times\cubeTk)\cap X, 
      \end{align}
     and 
        \begin{equation}
        \label{energytildevVSv}
        \limsup_{\e\to 0} (F_\e(w_\e) - F_\e(u_\e))\leq o(1)
        \end{equation}
   as $\delta \to 0$.
  \end{lemma}
   \begin{proof}
   Fixed $N\in\mathbb{N}$ and $\delta\in (0,1/4)$. For $h\in\{0,\cdots, N\}$, we set
       \begin{equation}
       \notag
       A_{h}:=\Bigl\{ x\in A\hspace{0.03cm}:\hspace{0.03cm} \dist(x, A(\delta))< h{\delta\over N} \Bigr\}.
       \end{equation}
    For $h\in\{0, \cdots, N-1 \}$, let $\phi^h_d$ be a cut-off function between $A_h$ and $A_{h+1}$ with $|\nabla\phi^h_d|\leq2N/\delta$ and let $w_\e$ be a function defined by
        \begin{equation}
        \label{deftildeve}
        w^\e_i = w_\e (\e i^d, \e i^k) := \phi^h_d(\e i^d)u^\e_i + (1-\phi^h_d(\e i^d))v^\e_i.
        \end{equation}
     Since both $u_\e$ and $v_\e$ converge to $u$ with respect to convergence given by \eqref{defconvergence}, we also deduce that $w_\e$ converges to $u$ with respect to \eqref{defconvergence}. By adding and subtracting the term $\phi^h_d(\e i^d)u^\e_j+(1-\phi^h_d(\e i^d))v^\e_j$, we get that 
          \begin{align}
          w^\e_i  -w^\e_j &= \phi^h_d(\e i^d) (u^\e_i- u^\e_j)+ (1-\phi^h_d(\e i^d))(v^\e_i-v^\e_i)+ (\phi^h_d(\e i^d)-\phi^h_d(\e j^d))(u^\e_j- v^\e_j). \label{difftildeve-ij}
          \end{align}
     For $h\in\{1,\dots, N-2\}$, we set 
          \begin{equation*}
          S_h^{d} := A_{h+1}\setminus A_h,
          \end{equation*}
      so that $A= A_h\cup A\setminus\overline{A}_{h+1}\cup S_h^d$. In order to estimate the energy
          \begin{equation}
          \notag
          \sum_{i,j\in ({1\over \e}A\times\cubeTk)\cap X}\e^{d-2}a_{ij}(w^\e_i-w^\e_j)^2,
          \end{equation}
      we separately evaluate the following cases
         \begin{itemize}
         \item[i)] $i,j\in ({1\over \e}A_h\times\cubeTk)\cap X$;
         \item[ii)] $i,j\in ({1\over \e}(A\setminus\overline{A}_{h+1} )\times\cubeTk)\cap X$;
         \item[iii)] $i\in ({1\over \e}S_h^d\times\cubeTk)\cap X$ and  $j\in ({1\over \e}A\times\cubeTk)\cap X$;
         \item[iv)] $i\in ({1\over \e}(A_h\cup( A\setminus\overline{A}_{h+1} ))\times\cubeTk)\cap X$ and $j\in({1\over \e}S_h^d\times\cubeTk)\cap X
          $;
         \item[v)] $i\in ( {1\over \e}A_h\times\cubeTk )\cap X$ and $j\in ({1\over \e}(A\setminus \overline{A}_{h+1} )\times\cubeTk)\cap X$;
         \item[vi)]  $i\in ( {1\over \e}(A\setminus \overline{A}_{h+1} )\times\cubeTk )\cap X$ and $j\in({1\over \e} A_h\times\cubeTk )\cap X$
         \end{itemize}
as follows

i) In view of definition \eqref{deftildeve}, we deduce that 
            \begin{align}
            \sum_{i,j\in ({1\over \e}A_h\times\cubeTk)\cap X}\e^{d-2}a_{ij}(w^\e_i-w^\e_j )^2 &= \sum_{i,j\in ({1\over \e}A_h\times\cubeTk)\cap X}\e^{d-2}a_{ij}(u^\e_i-u^\e_j )^2\notag\\
            &\leq \sum_{i,j\in ({1\over \e}A\times\cubeTk)\cap X}\e^{d-2}a_{ij}(u^\e_i-u^\e_j )^2.\label{intQhQh}
            \end{align}
            \smallskip
            
        ii)  We have that 
               \begin{align}
              \sum_{i,j\in ({1\over \e}(A\setminus \overline{A}_{h+1} )\times\cubeTk)\cap X}\e^{d-2}a_{ij}(w^\e_i-w^\e_j )^2& =\sum_{i,j\in ({1\over \e}(A\setminus \overline{A}_{h+1} )\times\cubeTk)\cap X}\e^{d-2}a_{ij} (v^\e_i -v^\e_j)^2\notag\\
              & \leq \sum_{i,j\in ({1\over \e}(A\setminus \overline{A(\delta)} )\times\cubeTk)\cap X }\e^{d}a_{ij} (v^\e_i-v^\e_j)^2.\label{f2}
               \end{align}
          %where we have set $A_\delta:=\{x^d\in Q_d\hspace{0.03cm}:\hspace{0.03cm} \dist(x^d, \partial Q_d)<\delta \}$.
          In view of definition of $v^\e_i$ given by \eqref{def:ve} and since $\e i^d + \e Q_d = \e(i^d-j^d) + \e j^d + \e Q_d$, we deduce that 
              \begin{align}
              |v^\e_i-v^\e_j|^2& = \left|\fint _{\e i^d+\e Q_d} u(x)dx - \fint _{\e j^d+\e Q_d} u(x)dx\right|^2\notag\\
              &= \left|\fint _{\e j^d+\e Q_d} (u(x+\e(i^d-j^d)) - u(x) )dx\right|^2\label{estveivej}.
              \end{align}
           Since $u\in H^1(\Omega)$, we have that 
               \begin{align}
               u(x+\e(i^d-j^d)) - u(x)  &= \int_{0}^{1} {\partial u\over \partial t} (x+\e t(i^d-j^d))dt\notag\\
               & = \int_{0}^{1} \nabla u (x+\e t(i^d-j^d)) \cdot \e (i^d-j^d)dt.\notag
               \end{align}
           This, combined with \eqref{estveivej} and the Fubini theorem,  implies that
               \begin{align}
               |v^\e_i-v^\e_j|^2 &=\left|\fint _{\e j^d+\e Q_d}\int_{0}^{1} \nabla u (x+\e t(i^d-j^d)) \cdot \e (i^d-j^d)dtdx \right|\notag\\
               &= {1\over \e^d }\left|\int_{0}^{1}\int _{\e j^d+\e Q_d} \nabla u (x+\e t(i^d-j^d)) \cdot \e (i^d-j^d)dxdt \right|\notag\\
               &\leq \e^{2-d} |i^d-j^d|^2 \int_{0}^{1}\int _{\e j^d+\e Q_d} |\nabla u (x+\e t(i^d-j^d))|^2 dxdt\notag\\
               &\leq \e^{2-d} T^2\int_{0}^{1}\int _{\e i^d+\e Q_d} |\nabla u (x))|^2 dxdt,\label{estveivej1}
               \end{align}
            where in the last inequality we have used the fact that the nodes $i$ and $j$ interact at most at distance $T$.
             In view of the assumption of finite range along with estimate above, from \eqref{f2}, it follows that
              \begin{align}
              \sum_{i,j\in({1\over \e} (A\setminus \overline{A}_{h+1}) \times\cubeTk)\cap X }\e^{d-2}a_{ij} (w^\e_i-w^\e_j)^2 &\leq C \sum_{i\in({1\over \e}(A\setminus \overline{A(\delta)}) \times\cubeTk)\cap X}\int_{\e i^d+\e Q_d}|\nabla u(x)|^2dx\notag\\
              &\leq C \int_{A\setminus A(2\delta)}|\nabla u(x)|^2dx,\label{intQd-h+1}
              \end{align}
        where the constant $C$ is due to the fact that a fixed node $i\in {1\over \e}(A\setminus \overline{A(\delta)}) \times\cubeTk)\cap X$ interacts with a finite number of nodes $j\in {1\over \e}(A\setminus \overline{A(\delta)}) \times\cubeTk)\cap X$.
 \smallskip
 
            iii) First note that due to the assumption of finite range, if $\e i^d\in S_h^d$, then $\e j^d\in \widehat{S}^d_h:=S_{h-1}^d\cup S_h^d\cup S_{h+1}^d$.  This  combined with  \eqref{difftildeve-ij} and the Jensen inequality implies that  
            \begin{align}
            &\sum_{\substack{i\in ({1\over \e}S_h^d\times\cubeTk)\cap X\\ j\in ({1\over \e}A\times\cubeTk)\cap X}}\e^{d-2}a_{ij}(w^\e_i-w^\e_j )^2\notag = \sum_{\substack{i\in ({1\over \e}S_h^d\times\cubeTk)\cap X\\ j\in ({1\over \e}\widehat{S}_h^d\times\cubeTk)\cap X}}\e^{d-2}a_{ij}(w^\e_i-w^\e_j )^2\notag\\
            &\quad \leq C \sum_{\substack{i\in ({1\over \e}S_h^d\times\cubeTk)\cap X\\ j\in ({1\over \e}\widehat{S}_h^d\times\cubeTk)\cap X}}\e^{d-2}a_{ij}(u^\e_i-u^\e_j )^2 + C\sum_{\substack{i\in ({1\over \e}S_h^d\times\cubeTk)\cap X\\ j\in ({1\over \e}\widehat{S}_h^d\times\cubeTk)\cap X}}\e^{d-2}a_{ij}(v^\e_i-v^\e_j )^2\notag\\
            &\quad\quad + C\sum_{\substack{i\in ({1\over \e}S_h^d\times\cubeTk)\cap X\\ j\in ({1\over \e}\widehat{S}_h^d\times\cubeTk)\cap X}}\e^{d-2}a_{ij}(\phi^h_d(\e i^d)-\phi^h_d(\e j^d))^2(u^\e_j-v^\e_j )^2.\label{intShd_1}
            \end{align} 
        Due to the fact that $|\nabla\phi^h_d|\leq 2N/\delta$, the last integral in \eqref{intShd_1} can be estimated as follows
           \begin{align}
           & \sum_{\substack{i\in ({1\over \e}S_h^d\times\cubeTk)\cap X\\ j\in ({1\over \e}\widehat{S}_h^d\times\cubeTk)\cap X}}\e^{d-2}a_{ij}(\phi^h_d(\e i^d)-\phi^h_d(\e j^d))^2(u^\e_j-v^\e_j )^2\notag\\
           &\quad\leq {N^2\over \delta^2} \sum_{\substack{i\in ({1\over \e}\widehat{S}_h^d\times\cubeTk)\cap X\\ j\in ({1\over \e}\widehat{S}_h^d\times\cubeTk)\cap X}}\e^{d-2}a_{ij}(u^\e_j-v^\e_j )^2\leq 3{N^2\over \delta^2}\sum_{i,j\in ({1\over \e}A\times\cubeTk)\cap X}\e^{d}a_{ij}(u^\e_j-v^\e_j )^2. \notag%\label{intShd_2}
           \end{align}
          In order to estimate the first two integrals in \eqref{intShd_1}, we may choose $h\in \{1,\dots,N-2 \}$  such that   
             \begin{align}
             &\sum_{\substack{i\in ({1\over \e}S_h^d\times\cubeTk)\cap X\\ j\in ({1\over\e}\widehat{S}_h^d\times\cubeTk)\cap X}}\e^{d-2}a_{ij}\left[(u^\e_i-u^\e_j )^2+(v^\e_i-v^\e_j )^2\right]\notag\\
             &\hspace{0.5cm}\leq {1\over N-2} \sum_{i,j\in ({1\over \e}A\times\cubeTk)\cap X}\e^{d-2}a_{ij}(u^\e_i-u^\e_j )^2+{1\over N-2} \sum_{i,j\in ({1\over \e}A\times\cubeTk)\cap X}\e^{d-2}a_{ij}(v^\e_i-v^\e_j )^2\notag\\
             &\hspace{0.5cm}\leq {1\over N-2} \sum_{i,j\in ({1\over \e}A\times\cubeTk)\cap X}\e^{d-2}a_{ij}(u^\e_i-u^\e_j )^2+{C\over N-2}\int_{A\setminus A(2\delta)} |\nabla u(x)|^2dx,\notag
              %&\hspace{0.5cm} \leq {1\over N-2} \sum_{i,j\in ({1\over \e}A\times\cubeTk)\cap X}\e^{d-2}a_{ij}(u^\e_i-u^\e_j )^2 +{C\over N-2}\fint_{\e i^d+\e Q^{x_0}_d} |\nabla u(x)|^2dx,\notag %\label{intShd_3}
             \end{align}  
         where  we have used \eqref{estveivej1} and the assumption of finite range.
         This, combined with \eqref{intShd_1}, leads us to 
             \begin{align}
             &\sum_{\substack{i\in ({1\over \e}S_h^d\times\cubeTk)\cap X\\ j\in ({1\over \e}A\times\cubeTk)\cap X}}\e^{d-2}a_{ij}(w^\e_i-w^\e_j )^2\leq C  {1\over N-2} \sum_{i,j\in ({1\over \e}A\times\cubeTk)\cap X}\e^{d-2}a_{ij}(u^\e_i-u^\e_j )^2\notag\\
             &\hspace{1.2cm}+{C\over N-2}\int_{A\setminus A(2\delta)} |\nabla u(x)|^2dx+ C{N^2\over \delta^2}\sum_{\substack{i\in ({1\over \e}S_h^d\times\cubeTk)\cap X\\ j\in ({1\over \e}A\times\cubeTk)\cap X}}\e^{d}a_{ij}(u^e_j-v^\e_j )^2.\label{intShd}
             \end{align}
              \smallskip
            
         iv) Note that 
             \begin{align}
             &\sum_{\substack{i\in ({1\over \e}( A_h\cup A\setminus\overline{A}_{h+1} )\times\cubeTk)\cap X\\ j\in ({1\over \e}S_h^d\times\cubeTk)\cap X}}\e^{d-2}a_{ij}(w^\e_i-w^\e_j )^2\leq\sum_{\substack{i\in ({1\over \e} A\times\cubeTk)\cap X\\ j\in ({1\over \e}S_h^d\times\cubeTk)\cap X}}\e^{d-2}a_{ij}(w^\e_i-w^\e_j)^2,\notag  
             \end{align}
        so that, the same argument as for iii) can be performed, obtaining estimate \eqref{intShd}.
        
\begin{comment}
    \item[v)] We get that 
            \begin{align}
            &\sum_{\substack{(\eta_\e i^d, i^k)\in ((Q_d)^h\times\cubeTk)\cap X_{\eta_\e}\\ (\eta_\e j^d, j^k)\in ((Q_d\setminus\overline{(Q_d)}^{h+1}  )\times\cubeTk)\cap X_{\eta_\e}}}\eta_\e^{d-2}a_{ij}(\tilde{v}_\e(\eta_\e i^d, \e i^k)-\tilde{v}_\e(\eta_\e j^d, \e j^k) )^2\notag\\
            &\quad =\sum_{\substack{(\eta_\e i^d, i^k)\in ((Q_d)^h\times\cubeTk)\cap X_{\eta_\e}\\ (\eta_\e j^d, j^k)\in ((Q_d\setminus\overline{(Q_d)}^{h+1}  )\times\cubeTk)\cap X_{\eta_\e}}}\eta_\e^{d-2}a_{ij}(v_\e(\eta_\e i^d, \e i^k)-\nabla u(x_0)\cdot\eta_\e j^d )^2.\label{intQh_Qd-h+1}
            \end{align}
        \item[vi)] We can repeat the same arguments as in the previous case.
\end{comment}     
       
\smallskip
               In view of the finite-range assumption, the points belonging to sets of items (v) and (vi) do not have any interaction since $\delta/N>>\e T$.

\begin{comment}
      Now, set
         \begin{align}
         G^h_\e&:=\sum_{\substack{(\eta_\e i^d, i^k)\in (S^h_d\times\cubeTk)\cap X_{\eta_\e}\\ (\eta_\e j^d, j^k)\in (Q_d\times\cubeTk)\cap X_{\eta_\e}}}\eta_\e^{d-2}a_{ij}(v_\e(\eta_\e i^d, \e i^k)-v_\e(\eta_\e j^d, \e j^k) )^2\notag\\
         &\quad +\sum_{\substack{(\eta_\e i^d, i^k)\in (S^h_d\times\cubeTk)\cap X_{\eta_\e}\\ (\eta_\e j^d, j^k)\in (Q_d\times\cubeTk)\cap X_{\eta_\e}}}\eta_\e^{d-2}a_{ij}(\nabla u(x_0)\cdot\eta_\e i^d-\nabla u(x_0)\cdot\eta_\e j^d )^2\notag\\
          &\quad +\sum_{\substack{(\eta_\e i^d, i^k)\in ((Q_d)^h\times\cubeTk)\cap X_{\eta_\e}\\ (\eta_\e j^d, j^k)\in ((Q_d\setminus\overline{(Q_d)}^{h+1}  )\times\cubeTk)\cap X_{\eta_\e}}}\eta_\e^{d-2}a_{ij}(v_\e(\eta_\e i^d, \e i^k)-\nabla u(x_0)\cdot\eta_\e j^d )^2.\notag
         \end{align}
\end{comment} 
   
         Gathering estimates \eqref{intQhQh}, \eqref{intQd-h+1} and \eqref{intShd}, we obtain that, for $h\in\{1,\dots, N-2\}$, 
          \begin{align}
          &\sum_{i,j\in ({1\over \e} A\times\cubeTk)\cap X}\e^{d-2}a_{ij} (w^\e_i-w^\e_j )^2\leq\sum_{i,j\in ({1\over \e}A\times\cubeTk)\cap X}\e^{d-2}a_{ij}(u^\e_i-u^\e_j )^2\notag\\
          &\hspace{2cm}+C\int_{A\setminus A(2\delta)} |\nabla u(x)|^2dx+ C  {1\over N-2} \sum_{i,j\in ({1\over \e}A\times\cubeTk)\cap X}\e^{d-2}a_{ij}(u^\e_i-u^\e_j )^2\notag\\
          &\hspace{2cm}+{C\over N-2}\int_{A\setminus A(2\delta)} |\nabla u(x)|^2dx+ C{N^2\over \delta^2}\sum_{\substack{i\in ({1\over \e}S_h^d\times\cubeTk)\cap X\\ j\in ({1\over \e}A\times\cubeTk)\cap X}}\e^{d}a_{ij}(u^\e_j-v^\e_j)^2. \label{estenergy}
          \end{align}
   
\begin{comment}
     We can choose an index $h\in\{1,\dots, N-1\}$ such that 
         \begin{align}
         G^h_\e&\leq {1\over N-1} \sum_{(\eta_\e i^d, i^k),(\eta_\e j^d, j^k)\in (Q_d\times\cubeTk)\cap X_{\eta_\e}}\eta_\e^{d-2}a_{ij}(v_\e(\eta_\e i^d, \e i^k)-v_\e(\eta_\e j^d, \e j^k) )^2\notag\\
          &\quad+{1\over N-1} \sum_{(\eta_\e i^d, i^k),(\eta_\e j^d, j^k)\in (Q_d\times\cubeTk)\cap X_{\eta_\e}}\eta_\e^{d-2}a_{ij}(\nabla u(x_0)\cdot\eta_\e i^d-\nabla u(x_0)\cdot\eta_\e j^d )^2\notag\\
          &\quad {\color{red}+ {1\over N-1} \sum_{(\eta_\e i^d, i^k) (\eta_\e j^d, j^k)\in (Q_d\times\cubeTk)\cap X_{\eta_\e}}\eta_\e^{d-2}a_{ij}(v_\e(\eta_\e i^d, \e i^k)-\nabla u(x_0)\cdot\eta_\e j^d )^2}\notag\\
          &\leq {C\over N-1}.\label{estGhe} 
         \end{align}
      Taking into account \eqref{estGhe}, we take the limit as $\e\to 0$ in \eqref{estenergy} obtaining
          \begin{align}
          &\limsup_{\e\to 0}
         \sum_{(\eta_\e i^d, i^k), (\eta_\e j^d, j^k)\in (Q_d\times\cubeTk)\cap X_{\eta_\e}}\eta_\e^{d-2}a_{ij}(\tilde{v}_\e(\eta_\e i^d, \e i^k)-\tilde{v}_\e(\eta_\e j^d, \e j^k) )^2\notag\\
         &\leq\limsup_{\e\to 0}  \sum_{(\eta_\e i^d, i^k), (\eta_\e j^d, j^k)\in (Q_d\times\cubeTk)\cap X_{\eta_\e}}\eta_\e^{d-2}a_{ij}(v_\e(\eta_\e i^d, \e i^k)-v_\e(\eta_\e j^d, \e j^k) )^2\notag\\
         &\quad +{\color{red}C\delta+{C\over N-1}},\notag
          \end{align}
\end{comment}
   Note that the last sum  vanishes as $\e \to 0$ since both $u_\e$ and $v_\e$ converge to $u$ with respect to convergence \eqref{defconvergence}. Hence, taking the limit as  $\e\to 0$ of \eqref{estenergy}, we obtain that 
       \begin{align*}
       \limsup_{\e\to 0} F_\e(w_\e) - F_\e(u_\e)&\leq C\int_{A\setminus A(2\delta)} |\nabla u(x)|^2dx + {C\over N-2} \liminf_{\e\to 0} F_\e(u_\e)\notag\\
       &\quad + {C\over N-2}\int_{A\setminus A(2\delta)} |\nabla u(x)|^2dx .\notag
       \end{align*}
    Letting first $N\to\infty$ and then  $\delta\to 0$, we get inequality \eqref{energytildevVSv} as desired.
   \end{proof}

\section{Homogenization}\label{Homsec}
This section is devoted to the limit analysis as $\e\to 0$ of the family of functionals $F_\e: {\cal C}_\e(\Omega)\to [0, \infty)$ defined as
      \begin{equation}
      \label{extendedenergy}
      F_\e(u) :=      
      \sum_{i,j\in ({1\over\e}\Omega\times \cubeTk)\cap X}\e^{d-2}a_{ij}(u_i-u_j)^2,
      \end{equation}
where $\Omega$ is a bounded open set of $\R^d$ with Lipschitz boundary and ${\cal C}_\e(\Omega)$ is given by \eqref{defCeOmegaX}. This is done through the computation of the corresponding $\Gamma$-limit with respect to convergence \eqref{convsuX}.

   \begin{thm}
   \label{thm:homogenization}
   The family of functionals \eqref{extendedenergy} $\Gamma$-converges to $F_\hom: H^1(\Omega)\to\R$ defined by 
         \begin{equation}
         F_\hom(u):=\int_{\Omega}A_\hom\nabla u\cdot\nabla udx,\notag
         \end{equation}
    where
       \begin{eqnarray}
       \label{periodicform}\nonumber
       A_\hom z\cdot z&:=&{1\over T^d}\min\biggl\{\sum_{i\in(\cubeTd\times\cubeTk)\cap X}\sum_{j\in( \R^d\times\cubeTk)\cap X}a_{ij}(u_i-u_j)^2\hspace{0.03cm}:\hspace{0.03cm}\\
       &&\hskip3cm u_i- z\cdot i^d\quad  \mbox{\rm is } T\mbox{\rm-periodic in }e_1,\dots,e_d \biggr\}.
       \end{eqnarray}
   \end{thm}
In this formula we interpret $u_i=z\cdot i^d$ as the discrete interpolation of the affine function $z\cdot x^d$, with $x^d\in\R^d$.

\subsection{Proof of the lower bound}
We prove the  lower-bound inequality for the family $F_\e$ using the {\it blow-up method} introduced by Fonseca and M\"{u}ller \cite{FM} (see also \cite{BMS}). 

Let $u_\e$ be a sequence with equi-bounded energy $F_\e(u_\e)$ and such that $u_\e$ converge to $u\in H^1(\Omega)$. Let the sequence of positive  measures $\lambda_\e$ be defined as
    \begin{equation}
    \notag
    \lambda_\e:= \sum_{i\in ({1\over \e}\Omega\times \cubeTk)\cap X}\biggl(\sum_{j\in({1\over \e}\Omega\times \cubeTk)\cap X}\e^{d-2}a_{ij}(u^\e_i-u^\e_j)^2  \biggr)\delta_{\e i},
    \end{equation} 
where $\delta_x$ is the Dirac measure concentrated at $x$. The $d$-dimensional measure $\mu_\e$ 
 is defined by
    \begin{equation}
    \notag
    \mu_\e(B) :=\lambda_\e (B\times\e\cubeTk) = \sum_{i\in ({1\over \e}B\times \cubeTk)\cap X}\sum_{j\in({1\over \e}\Omega\times \cubeTk)\cap  X}\e^{d-2}a_{ij}(u^\e_i-u^\e_j)^2
    \end{equation}for Borel sets $B$ of $\R^{d}$.
Note that $\mu_\e(B)$ takes into account  interactions between the nodes with projection in $B$ and the ones in all $X$, but, in view of the equi-boundedness of ${\cal E}$, which is a finite-range assumption, we can limit the interactions between the nodes with projection in $B$ and those with projection in an $\e R$-neighbourhood of $B$. \par
Since $\mu_\e(\Omega)=F_\e(u_\e)$ and thanks to the equi-boundedness of $F_\e(u_\e)$, the measures $\mu_\e$ are also equi-bounded, so that, up to subsequences, we deduce that 
    \begin{equation*}
    \mu_\e\weakstar\mu,
    \end{equation*} 
where $\mu$ is a $d$-dimensional positive measure on $\Omega$. The Radon-Nikodym decomposition of the limit measure $\mu$ with respect to the $d$-dimensional Lebesgue measure ${\cal L}^d$ enables us to write that 
      \begin{equation}
      \notag
      \mu ={d\mu\over dx}{\cal L}^d+\mu^{s},
      \end{equation}
with $\mu^s\perp{\cal L}^d$.  Note that the positiveness of $\mu$ ensures that its singular part $\mu^s$ is positive as well. \par Now, we perform a local analysis. Let $x_0\in\Omega$ be a Lebesgue point for $\mu$ with respect to ${\cal L}^d$; {\it i.e.},
     \begin{equation}
     \label{lebesguepoint}
     {d\mu\over dx}(x_0)=\lim_{\rho\to 0}{\mu( Q_{\rho,d}(x_0))\over {\cal L}^d( Q_{\rho,d}(x_0)) } =\lim_{\rho\to 0} {\mu( Q_{\rho,d}(x_0))\over \rho^d},
     \end{equation}  
with $Q_{\rho,d}(x_0):= x_0+[0,\rho)^d$.
Thanks to the Besicovitch Derivation Theorem, ${\cal L}^d$-almost every $x_0\in\Omega$ is a Lebesgue point for $\mu$ with respect to ${\cal L}^d$. Moreover, in view of \cite[Theorem 3.4.2]{Z89}, we have that, up to a set of zero Lebesgue measure,  $x_0$ is a point such that 
     \begin{equation}
     \label{propLebpoint}
     \lim_{\rho\to 0} {1\over \rho}\biggl({1\over\rho^d}\int_{ Q_{\rho,d}(x_0)}|u(x)-u(x_0)- \nabla u(x_0)\cdot (x-x_0)|^2 dx  \biggr)^{1/2}=0.
     \end{equation}
In other words, performing the change of variables $x=\rho y +x_0$ in the above integral, we have that
    \begin{equation*}
    {u(\rho y+x_0) - u(x_0)\over \rho} \to \nabla u(x_0)\cdot y\quad\mbox{in } L^2(Q_d).
    \end{equation*}
 For all $\rho\to 0$ but a countable set, we have that $\mu( Q_{\rho,d}(x_0))=0$ and hence for such $\rho$ we have that 
   \begin{equation}
   \label{limitmeasure}
   \mu(Q_d(x_0))=\lim_{\e\to 0} \mu_\e( Q_{\rho,d}(x_0)).
   \end{equation}
Therefore, from \eqref{lebesguepoint}, it follows that
   \begin{equation}
   \notag
   {d\mu\over dx}(x_0)=\lim_{\rho\to 0}\lim_{\e\to 0} { \mu_\e(Q_{\rho,d}(x_0))\over \rho^d  }.
   \end{equation}
Now, we perform the blow-up argument.  Since $x_0\in\Omega$ is a Lebesgue point and due to a  diagonalization argument on \eqref{lebesguepoint} and \eqref{limitmeasure}, there exists a sequence $\rho_\e\to 0$ as $\e\to 0$ such that $\rho_\e >\!>\e$ and the following equalities
    \begin{equation}
    \label{lebpoint1}
    {d\mu\over dx}(x_0)=\lim_{\e\to 0} { \mu_\e( Q_{\rho_\e,d}(x_0))\over \rho_\e^d  },
    \end{equation}
and 
    \begin{equation}
    \label{lebpoint2}
    \lim_{\e\to 0}{1\over \rho_\e} \int_{Q_{\rho_\e,d}(x_0)\times\cubeTk} |u_\e-u|(x)\chi_{\cup_{i\in X}\e Q_{d}^{i^d}\times Q^{i^k}_k}dx=0,
    \end{equation}
hold. Thanks to the link between the measure $\mu_\e$ and the energy $F_\e$, equality \eqref{lebpoint1} can be re-written as
   \begin{equation}
   \notag
   {d\mu\over dx}(x_0)=\lim_{\e\to 0} { 1\over \rho_\e^d  }\sum_{i\in( Q_{{\rho_\e\over \e},d}\left({x_0 \over \e} \right)\times \cubeTk)\cap X}\sum_{j\in({1\over \e}\Omega\times\cubeTk)\cap X}\e^{d-2}a_{ij}(u^\e_i-u^\e_j)^2. 
   \end{equation}
Now, the aim is to estimate the limit above. First, note that since the coefficients $a_{ij}$ are positive, we can consider only interactions taking place between nodes inside the cube $Q_{{\rho_\e\over \e},d}(x_0/\e)\times \cubeTk$, so that 
    \begin{equation}
    \label{e1}
    { \mu_\e( Q_{\rho_\e,d}(x_0))\over \rho_\e^d  } \geq{ 1\over \rho_\e^d  }\sum_{i,j\in( Q_{{\rho_\e\over \e},d}({x_0\over \e})\times\cubeTk)\cap X}\e^{d-2}a_{ij}(u^\e_i-u^\e_j)^2.
    \end{equation}
We need to modify $u_\e$ in order to define a function $v_\e$  converging to the affine function $\nabla u(x_0)\cdot x^d$ in $L^2(Q_d)$. %To this end, we re-scale  $Q_{{\rho_\e\over \e},d}(x_0/ \e)\times\cubeTk\cap X$ to $Q_d\times\cubeTk\cap X_{\eta_\e}$, where $X_{\eta_\e}$ is the set $X$ rescaled to  $\eta_\e\Z^d\times \{0,\dots, T-1\}^k$ with $\eta_\e:=\e/\rho_\e $.
To that end, let $\eta_\e={\e\over\rho_\e}$, and let $X_{\eta_\e, \e}$ be the set $X$ rescaled to $\eta_\e\Z^d\times \e \{0,\dots, T-1\}^k$.  We define  $v^\rho_\e$ on $(Q_d\times\cubeTk)\cap X_{\eta_\e, \e}$ by
     \begin{equation}
     \label{defvrhoe}
     v^\rho_\e(\eta_\e i^d, \e i^k):= {u_\e(\e i^d+ x_0, \e i^k) - u(x_0) \over \rho_\e},
     \end{equation}
where $u_\e$ is defined on $(Q_{\rho_\e, d}(x_0)\times\cubeTk)\cap \e X$.
Note that since $u_\e$ is a function in ${\cal C}_\e(\Omega)$, $v^\rho_\e$ can be identified  with a piecewise-constant function on $\eta_\e Q_d^{i^d}\times \e Q_k^{i^k}$ if $(i^d, i^k)\in X$. For $x=(x^d,x^k)\in\R^d\times \cubeTk$, we set $w_0(x):= \nabla u(x_0)\cdot x^d$ and we show that 
   \begin{equation}
   \label{limvrhoe-aff}
   \lim_{\e \to 0} \int_{Q_d\times \cubeTk} |v^\rho_\e(x)-w_0(x)|^2\chi_{\cup_{i\in X}\eta_\e Q_d^{i^d}\times Q_k^{i^k}}(x)dx=0.
   \end{equation}
To this end, we introduce the function $u_0$ given by $u_0(x^d):= u(x_0) + w_0(x) $. Hence, 
   \begin{align}
   \notag
   u(x_0) = u_0(\rho x^d) - \rho\nabla u(x_0)\cdot x^d = u_0(\rho x^d) - \rho w_0(x).
   \end{align} 
This, combined with \eqref{defvrhoe}, implies that
   \begin{align}
  &\int_{Q_d\times \cubeTk} |v^\rho_\e(x)-w_0(x)|^2\chi_{\cup_{i\in X}\eta_\e Q_d^{i^d}\times Q_k^{i^k}}(x)dx\notag\\
  &\hspace{1cm}=\int_{Q_d\times \cubeTk}\biggl| {u_\e (\rho_\e x^d+x_0, \e x^k) - u(x_0)\over \rho_\e}-w_0(x)\biggr|^2\chi_{\cup_{i\in X}\eta_\e Q_d^{i^d}\times Q_k^{i^k}}(x)dx\notag\\
  &\hspace{1cm}=\int_{Q_d\times \cubeTk}\biggl| {u_\e (\rho_\e x^d+x_0, \e x^k) - u_0(\rho_\e x^d)\over \rho_\e}\biggr|^2\chi_{\cup_{i\in X}\eta_\e Q_d^{i^d}\times Q_k^{i^k}}(x)dx\notag\\
  &\hspace{1cm}\leq\int_{Q_d\times \cubeTk}\biggl| {u_\e (\rho_\e x^d+x_0, \e x^k) -u(\rho_\e x^d+x_0)\over \rho_\e}\biggr|^2\chi_{\cup_{i\in X}\eta_\e Q_d^{i^d}\times Q_k^{i^k}}(x)dx\notag\\
  &\hspace{1.5cm}+\int_{Q_d\times \cubeTk}\biggl| {u (\rho_\e x^d+x_0) -u_0(\rho_\e x^d)\over \rho_\e}\biggr|^2\chi_{\cup_{i\in X}\eta_\e Q_d^{i^d}\times Q_k^{i^k}}(x)dx\label{e2}.
   \end{align}
The first integral in \eqref{e2} goes to $0$ as $\e \to 0$. Indeed, due to the change of variables $y^d=\rho_\e x^d +x_0$, we deduce that  
    \begin{align}
    \int_{Q_d\times \cubeTk}&\biggl| {u (\rho_\e x^d+x_0) -u_0(\rho_\e x^d)\over \rho_\e}\biggr|^2\chi_{\cup_{i\in X}\eta_\e Q_d^{i^d}\times Q_k^{i^k}}(x)dx^ddx^k\notag\\
    &\hspace{1cm}={1\over \rho_\e^{d+2}}\int_{Q_{\rho_\e, d}(x_0)\times \cubeTk}|u_\e (y^d, \e x^k)-u(y^d)|^2\chi_{\cup_{i\in X}\e Q_d^{i^d}\times Q_k^{i^k}}(x)dy^ddx^k,\notag 
    \end{align}
which vanishes as $\e\to 0$ thanks to \eqref{lebpoint2}. We  evaluate the second integral in \eqref{e2}.  Using again the change of variables $y^d=\rho_\e x^d+x_0$ and the definition of $u_0$, we have that 
    \begin{eqnarray*}
    &&\hskip-2cm\int_{Q_{\rho_\e, d}(x_0)\times \cubeTk}|u(y^d, \e x^k)-u(y^d)|^2\chi_{\cup_{i\in X}\e Q_d^{i^d}\times Q_k^{i^k}}(x)dy^ddx^k\\
    &\leq& T^k{1\over\rho_\e^2}{1\over \rho_\e^d}\int_{Q_{\rho_\e, d}(x_0)}| u (y^d) -u_0(y-x_0)|^2dy^d\\
    &= &T^k{1\over\rho_\e^2}{1\over \rho_\e^d}\int_{Q_{\rho_\e, d}(x_0)} |u(y^d)-u(x_0)-\nabla(x_0)\cdot(y-x_0)|^2dy^d.
    \end{eqnarray*}
Thanks to \eqref{propLebpoint}, it follows that also the  integral above vanishes as $\e\to 0$ so that we can conclude that \eqref{limvrhoe-aff} holds. Set
     \begin{equation*}
     v_\e(\eta_\e i^d, \e i^k):= v^{\rho_\e}_\e(\eta_\e i^d, \e i^k).
     \end{equation*}
Now, using Lemma \ref{lemma:treatmentBC}, we may modify the sequence $v_\e$ to get a new sequence  $\tilde{v}_\e$ which is equal to $\nabla u(x_0)\cdot \eta_\e i^d$ near the boundary $(\partial Q_d\times\cubeTk)\cap X_{\eta_\e}$, where $X_{\eta_\e}$ is the set $X$ rescaled to $\eta_\e\Z^d\times  \{0,\dots, T-1\}^k$, and
       \begin{align}
       &\limsup_{\e\to 0}\sum_{(\eta_\e i^d, i^k), (\eta_\e j^d, j^k)\in (Q_d\times\cubeTk)\cap X_{\eta_\e} }\eta^{d-2} a_{ij} (\tilde{v}_\e(\eta_\e i^d, \e i^k) - (\tilde{v}_\e(\eta_\e j^d, \e j^k))\notag\\
       &\hspace{1cm}\leq \limsup_{\e\to 0} \sum_{(\eta_\e i^d, i^k), (\eta_\e j^d, j^k)\in (Q_d\times\cubeTk)\cap X_{\eta_\e} }\eta^{d-2} a_{ij} (v_\e(\eta_\e i^d, \e i^k) - (v_\e(\eta_\e j^d, \e j^k))+o(1).\label{energytildevVSv}
       \end{align}

In order to simplify the notation, we may assume that $x_0\in\e T\Z^d$ so that we avoid the translation of the coefficients $a_{ij}$.
In view of \eqref{e1} and thanks estimate \eqref{energytildevVSv}, we have that
    \begin{align}
    {d\mu\over dx}(x_0)&\geq\limsup_{\e\to 0} \sum_{(\eta_\e i^d, i^k), (\eta_\e j^d, j^k)\in (Q_{d}({x_0 \over \rho_\e})\times \cubeTk)\cap X_{\eta_\e}}\eta_\e^{d-2}a_{ij}\left( {u^\e_i-u^\e_j\over \rho_\e}\right)^2 \notag\\
    &= \limsup_{\e\to 0} \sum_{(\eta_\e i^d, i^k), (\eta_\e j^d, j^k)\in (Q_{d}\times \cubeTk)\cap X_{\eta_\e}}\eta_\e^{d-2}a_{ij}(v_\e(\eta_\e i^d, \e i^k)- v_\e(\eta_\e j^d, \e j^k))^2\notag\\
    &\geq  \limsup_{\e\to 0} \sum_{(\eta_\e i^d, i^k), (\eta_\e j^d, j^k)\in (Q_{d}\times \cubeTk)\cap X_{\eta_\e}}\eta_\e^{d-2}a_{ij}(\tilde{v}_\e(\eta_\e i^d, \e i^k)- \tilde{v}_\e(\eta_\e j^d, \e j^k))^2\notag\\
    &\geq\limsup_{\e\to 0} \inf\Bigl\{\sum_{(\eta_\e i^d, i^k), (\eta_\e j^d, j^k)\in (Q_{d}\times \cubeTk)\cap X_{\eta_\e}}\eta_\e^{d-2}a_{ij}(w_\e(\eta_\e i^d, \e i^k)- w_\e(\eta_\e j^d, \e j^k))^2\hspace{0.03cm}:\notag\\
    &\hspace{2.5cm} w_\e(\eta_\e i^d, \e i^k)=\nabla u(x_0)\cdot \eta_\e i^d ,\quad \mbox{if }\dist(\eta_\e i^d, \partial Q_d )<2\eta_\e \sqrt{d} T\Bigr\}.\notag
    \end{align}  
Setting $K_\e=\lfloor 1/(\eta_\e T)\rfloor $, we have that 
   \begin{align}
   {d\mu\over dx}(x_0)&\geq\liminf_{\e\to 0}{1\over (K_\e T)^d} \inf\biggl\{\sum_{i, j\in (Q_{K_\e T,d}\times \cubeTk)\cap X}{1\over \eta_\e^2}a_{ij}(w_\e(\eta_\e i^d, \e i^k)- w_\e(\eta_\e j^d, \e j^k))^2\hspace{0.03cm}:  \notag\\
   &\hspace{3cm}\quad w_\e(\eta_\e i^d, \e i^k)=\nabla u(x_0)\cdot  \eta_\e i^d \mbox{ if }\dist( \eta_\e i^d, \partial Q_{K_\e T,d} )<2\eta_\e\sqrt{d}T  \biggr\}\notag\\
     &\geq\liminf_{\e\to 0}{1\over(K_\e T)^d} \inf\biggl\{\sum_{i, j\in (Q_{K_\e T,d}\times \cubeTk)\cap X}a_{ij}(\tilde{w}_i- \tilde{w}_i)^2\hspace{0.03cm}: \notag\\
     &\hskip5.5cm \tilde{w}_i=\nabla u(x_0)\cdot  i^d\mbox{ if }\dist( i^d, \partial Q_{k_\e T,d} )<2\sqrt{d}T \biggr\}\notag\\
            &=f_0 (\nabla u (x_0)), \notag
   \end{align}
   where we have set $\tilde{w}_i:=  w_\e(\eta_\e i^d, \e i^k)/\eta_\e$.  Therefore, for ${\cal L}^d$-almost every $x_0\in \Omega$, we have 
      \begin{align}
      {d\mu \over dx}(x_0)\geq A_\hom \nabla u(x_0)\cdot \nabla u(x_0).\notag
      \end{align}
  Integrating on $\Omega$, we conclude that
     \begin{equation}
     \notag
     \mu(\Omega)\geq \int_{\Omega}{d\mu(\Omega)\over dx}dx\geq \int_{\Omega}A_\hom \nabla u(x)\cdot \nabla u(x)dx.
     \end{equation} 
  Since $\mu_\e\weakstar\mu$, we have that $\liminf_{\e\to 0 }\mu_\e(\Omega)\geq\mu(\Omega)$. This implies that 
      \begin{equation}
      \notag
      \liminf_{\e\to 0} F_\e(u_\e)=\liminf_{\e\to 0}\mu_\e(\Omega)\geq \mu(\Omega) \geq \int_{\Omega}A_\hom \nabla u(x)\cdot \nabla u(x)dx = F_\hom (u),
      \end{equation}
    which concludes the proof of the lower bound.\par  It remains to prove that $f_0$ satisfies formula \eqref{periodicform}. First, we prove the existence of the limit.
  \begin{prop}
  There exists the limit
    \begin{align}
  f_0(z) = \lim_{K\to\infty}{1\over (KT)^d}\inf\biggl\{\sum_{i,j\in (Q_{KT,d}\times\cubeTk)\cap X}a_{ij}(u_i-u_j)^2\hspace{0.03cm}:\hspace{0.03cm}u_i = z\cdot i^d\quad  \mbox{ if }\dist(i^d, \partial Q_{KT,d})<2\sqrt{d}T  \biggr\},\label{ahf}
    \end{align}
  for $z\in\R$.
  \end{prop}
  \begin{proof}
  For fixed $K\in\mathbb{N}$ and $z\in\R^d$, we set 
    \begin{align}
    f^K_0(z) :={1\over (KT)^d}\inf\biggl\{\sum_{i,j\in (Q_{KT,d}\times\cubeTk)\cap X}a_{ij}(u_i-u_j)^2\hspace{0.03cm}:\hspace{0.03cm}u_i = z\cdot i^d\quad  \mbox{ if }\dist(i^d, \partial Q_{KT,d})<2\sqrt{d}T  \biggr\}.\notag
    \end{align} 
  Let $u^K$ be a function such that 
     \begin{align}
     \notag
     {1\over (KT)^d} \sum_{i,j\in(Q_{KT,d}\times\cubeTk)\cap X}a_{ij} (u^K_i-u^k_j)^2\leq f_0^K(z)+{1\over K},
     \end{align}
  and $u^K_i=z\cdot i^d$, if $\dist(i^d, \partial Q_{KT,d})<2\sqrt{d} T$. For $H>K$, we introduce the set of indices ${\cal I}:=\{l\in\Z^d\hspace{0.03cm}:\hspace{0.03cm} 0\leq (
   K +1)l_m<H, m=1,\dots, d \}$. We define
       \begin{align*}
       u^H_i:=
       \begin{dcases}
       u^K(i^d-l, i^k)+z\cdot l, &\quad (i^d, i^k)\in Q_{KT,d}^{l}\times\cubeTk,\hspace{0.1cm} l\in{\cal I},\\
       z\cdot i^d,&\quad \mbox{otherwise}.
       \end{dcases}
       \end{align*}
      We have that  
       \begin{align}
       f^H_0(z)&\leq {1\over (HT)^d} \sum_{i,j\in(Q_{HT,d}\times\cubeTk)\cap X}a_{ij}(u^H_i-u^H_j)^2\notag\\
       &= {1\over (HT)^d}\sum_{i\in(\cup_{l\in{\cal I}}Q^l_{KT,d}\times\cubeTk)\cap X} \sum_{j\in (Q_{HT,d}\times\cubeTk)\cap X } a_{ij}(u^K_i-u^H_j)^2\notag\\
       &\quad+ {1\over (HT)^d}\sum_{i\in [(Q_{HT,d}\setminus\cup_{l\in{\cal I}}Q^l_{K,d})\times\cubeTk]\cap X}\sum_{j\in (Q_{HT,d}\times\cubeTk)\cap X }a_{ij} (z\cdot i^d-u^H_j)^2.\label{estAHhom}
       \end{align}
  Due to the assumption of finite range, the interactions between nodes in $(\cup_{l\in{\cal I}}Q^l_{KT,d}\times Q_{T,k})\cap X$ and $[(Q_{HT,d}\setminus\cup_{l\in{\cal I}}Q^l_{KT,d})\times\cubeTk]\cap X$ do not take place. This implies that the first sum in \eqref{estAHhom} may be estimated as
     \begin{align}
     &\hskip-1cm{1\over (HT)^d}\sum_{i\in(\cup_{l\in{\cal I}}Q^l_{KT,d}\times\cubeTk)\cap X} \sum_{j\in (Q_{HT,d}\times\cubeTk)\cap X }  a_{ij}(u^K_i-u^H_j)^2\notag\\
     &={1\over (HT)^d}\sum_{i,j\in(\cup_{l\in{\cal I}}Q^l_{KT,d}\times\cubeTk)\cap X} a_{ij}(u^K_i-u^K_j)^2\leq {K^d\over H^d}\sum_{l\in{\cal I}}(f^K_0(z) + K^{-1})\notag\\
     &\leq {K^d\over H^d}\left\lfloor {H\over K+1}\right\rfloor^d ( f^H_0 (z) +  K^{-1})\leq {K^d\over( K+1)^d}(f^K_0 (z)+ K^{-1}).\label{AKhom1}
     \end{align}
  Using the assumption of finite range, the second sum in \eqref{estAHhom} may be estimated as
     \begin{align}
     {1\over (HT)^d}&\sum_{i\in [(Q_{HT,d}\setminus\cup_{l\in{\cal I}}Q^l_{KT,d})\times\cubeTk]\cap X}\sum_{j\in (Q_{HT,d}\times\cubeTk)\cap X }a_{ij} (z\cdot i^d-u^H_j)^2\notag\\
      &\leq{1\over (HT)^d}\sum_{i\in [(Q_{HT,d}\setminus\cup_{l\in{\cal I}}Q^l_{KT,d})\times\cubeTk]\cap X}\left(\sum_{j\in (Q_{HT,d}\times\cubeTk)\cap X} a_{ij}|z\cdot i^d-z\cdot j^d|^2 \right)\notag\\
      &\leq{C\over (HT)^d}\sum_{i\in [(Q_{HT,d}\setminus\cup_{l\in{\cal I}}Q^l_{KT,d})\times\cubeTk]\cap X}\left(\sum_{j\in (Q_{HT,d}\times\cubeTk)\cap X} a_{ij}|i^d-j^d|^2 \right)\notag\\
      &\leq {C\over (HT)^d}\left\lfloor {H\over  K +1}\right\rfloor^d.\label{AKhom2}
     \end{align}
  Combining \eqref{AKhom1} and \eqref{AKhom2}, from \eqref{estAHhom} it follows that 
      \begin{align}
      f^H_0 (z)&\leq  {K^d\over( K+1)^d}(f^K_0 (z)+ K^{-1})+ {1\over T^d(K+1)^d}.\notag
      \end{align}
  Taking first the limsup as $H\to \infty$ and then the lower limit as $K\to\infty$, we obtain
       \begin{align}
       \notag
       \limsup_{H\to\infty} f^H_0 (z)\leq \liminf_{K\to\infty} f^K_0 (z),
       \end{align}
   which concludes the proof.
  \end{proof}
Since we deal with  convex energies, asymptotic homogenization formula \eqref{ahf} can be reduced to a single periodic minimization problem. 
   \begin{prop}
   \label{prop:perminpb}
   We have that $f_0(z)$ defined by \eqref{ahf} coincides with $ f_\hom (z)$ defined as
             \begin{equation}
             \notag%\label{periodicform1}
            f_\hom (z):={1\over T^d}\inf\biggl\{\sum_{i\in(\cubeTd\times\cubeTk)\cap X}\sum_{j\in( \R^d\times\cubeTk)\cap X}a_{ij}(u_i-u_j)^2\hspace{0.03cm}:\hspace{0.03cm}u_i- z\cdot i^d\quad  \mbox{\rm is } T\mbox{\rm-periodic in }e_1,\dots,e_d \biggr\},
             \end{equation}
   for $z\in\R^d$.
   \end{prop}
   \begin{proof}
   Fix $z\in\R^d$.  First, we prove that $f_0 (z)\leq f_\hom (z)$. To this end, for $\delta>0$, let $u^\#$ be a function satisfying 
         \begin{align}
         \notag
         {1\over T^d}\sum_{i\in(\cubeTd\times\cubeTk)\cap X}\sum_{j\in (\R^d\times\cubeTk)\cap X} a_{i,j}(u^\#_i-u^\#_j)^2\leq f_\hom (z) +\delta
         \end{align}
   and $u^\#_i-z\cdot i^d$ is $T$-periodic in $e_1, \dots, e_d$. We define $u^\e_i= u_\e(\e i):= \e u^\#(i)$. Note that $u^\e_i$ converges to $z\cdot x^d$ with respect to the convergence given by \eqref{defconvergence}. Set $I^d:=\{l\in\Z^d\hspace{0.03cm}:\hspace{0.03cm} \e lT+\e\cubeTd\cap\Omega\neq\emptyset \}$. In view of Theorem \ref{thm:homogenization} and the periodicity of $a_{ij}$, we deduce that
       \begin{align}
       |\Omega|f_0 (z)&\leq \liminf_{\e\to 0} F_\e(u_\e)\notag\\
       &\leq \limsup_{\e\to 0}F_\e(u_\e)\leq \limsup_{\e\to 0} \sum_{l\in I^d}\sum_{i,j\in (Q_{T,d}^{l}\times\cubeTk)\cap X} \e^d a_{ij} (u^\#_i-u^\#_j)^2\notag\\
       &\leq \limsup_{\e\to 0}\sum_{l\in I^d}\e^d\sum_{i\in (Q_{T,d}\times\cubeTk)\cap X} \sum_{j\in (\R^d\times\cubeTk)\cap X}a_{ij}(u^\#_i-u^\#_j)^2\notag\\
       &\leq |\Omega|(f_\hom (z)+\delta).\notag 
       \end{align}
   From the arbitrariness of $\delta$, it follows the conclusion. \par It remains to show that $f_\hom (z)\leq f_0 (z)$. Let $v$ be a function defined on $(Q_{ KT, d}\times \cubeTk)\cap X$ such that $v_i=z\cdot i^d$ if $\dist(i^d, \partial Q_{KT, d})\leq 2\sqrt{d}T$. We define a function $u$ on $(\cubeTd\times\cubeTk)\cap X$ by
      \begin{align}
      \notag 
      u(i):={1\over K^d}\sum_{l\in\{0,\dots, K-1\}^d} v(i^d+lT, i^k).
      \end{align}
   With the help of Jensen's inequality combined with the assumption of finite range and the periodicity of $a_{ij}$, we deduce that 
      \begin{align}
       f_\hom (z) &\leq {1\over T^d}\sum_{i\in (\cubeTd\times\cubeTk)\cap X}\sum_{j\in (\R^d\times\cubeTk)\cap X}a_{ij}(u_i-u_j)^2\notag\\
       &\leq {1\over (KT)^{d}}\sum_{i\in (\cubeTd\times\cubeTk)\cap X}\sum_{j\in (\R^d\times\cubeTk)\cap X}\sum_{l\in\{0,\dots, K-1\}^d} a_{ij}(v(i^d+lT, i^k)-v(j^d+lT, j^k))^2\notag\\
       &={1\over (KT)^{d}}\sum_{i,j\in (Q_{KT,d}\times\cubeTk)\cap X} a_{ij}(v_i-v_j)^2.\notag
      \end{align}
   Taking the infimum, we get 
       \begin{align}
       \notag
       f_\hom (z)  &\leq{1\over (KT)^d} \inf\biggl\{\sum_{i,j\in (Q_{KT,d}\times\cubeTk)\cap X} a_{ij}(v_i-v_j)^2\hspace{0.03cm}:\hspace{0.03cm} v_i=z\cdot i^d  \mbox{ if }\dist(i^d, \partial Q_{KT,d})<2\sqrt{d}T \biggr\}.
       \end{align} 
   Then, passing to the limit as $K\to\infty$, we have the desired inequality which concludes the proof. 
   \end{proof}

\subsection{Proof of the upper bound}
We are going to prove the $\Gamma$-$\limsup$ inequality. The proof is independent of the blow-up result and it relies on a standard density argument by piecewise-affine functions (see \cite[Remark 1.29]{Brai02}). We consider the case when the target function $u$ is piecewise-affine and we assume that the gradient of $u$  takes $\lambda$ values, for some $\lambda$ positive integer.  For fixed $z_1,\dots, z_\lambda\in\R$, we define 
    \begin{equation}
    \notag
    %S:=\{x^d\in\Omega\hspace{0.03cm}:\hspace{0.03cm} z_1\cdot x^d=z_2\cdot z^d \}, \qquad
    \Omega_q:=\{x^d\in\Omega\hspace{0.03cm}:\hspace{0.03cm} u(x^d)=z_q\cdot x^d +c_q\},
    \end{equation}
for $q=1,\dots, \lambda$ (with $c_q$ some constant).  

We fix one such $q$.
We choose $w^q\in{\cal C}_\e(Q_{T,d})$ such that $w^q_i-z\cdot i^d$ is $T$-periodic in $e_1,\dots, e_d$ and
    \begin{align}
    \sum_{i,j\in (Q_{T,d}\times\cubeTk)\cap X} a_{ij}(w^q_i-w^q_j)^2 =  f_\hom(z_q).\notag
    \end{align} 
For  any $q=1,\dots, \lambda$, we define  $u^{\e,q}_i:= u_\e^q(\e i) = \e w^q(i)+c_q$. In view of Lemma \ref{lemma:treatmentBC}, we may modify the sequence $u^{q}_\e$ to obtain a new sequence $v^{q,\delta}_\e$ converging to $z_q\cdot x^d +c_q$ with respect to convergence \eqref{defconvergence} such that
   \begin{align}
   v^{\e,q,\delta}_i&=z_q\cdot \e i^d +c_q, \quad \mbox{if }  i\in ((\Omega_q\setminus\overline{\Omega_q(\delta)})\times\cubeTk)\cap X,\notag\\
   v^{\e,q,\delta}_i&=u^{\e,q}_i, \hspace{1.55cm}\mbox{if }  i\in (\Omega_q(2\delta)\times\cubeTk)\cap X,\notag
   \end{align} 
%where $\Omega_q(\delta)$ is defined by
% \begin{equation}
%    \Omega_q(\delta) :=\{x^d\in\Omega_q\hspace{0.03cm}:\hspace{0.03cm} \dist (x^d, \partial\Omega_q)>\delta \},\notag
%    \end{equation}
%for some $\delta>0$. We also have that
and 
      \begin{equation}
      \limsup_{\e\to 0} F^q_\e(v^{q,\delta}_\e) \leq \limsup_{\e\to 0} F^q_\e(u^q_\e) + o(1)\label{estOmegaq}
      \end{equation}
     as  $\delta\to 0$, where $F^q_\e$ is the functional defined as in \eqref{extendedenergy} with $\Omega_q$ in the place of $\Omega$.
      
\par Now, we estimate $F^q_\e(u_\e^q)$. To that end, for $q=1,\dots, \lambda$, we introduce the set of indices $
    {\cal I}^d_{\e,q}:= \{l\in\Z^d\hspace{0.03cm}:\hspace{0.03cm} \e Q^{l}_{T,d} \cap \Omega_q \neq\emptyset\},$ and  we deduce that 
 %Since $\Omega_q\subset\cup_{l\in{\cal I}^d_{\e,q}}\e Q^{lT}_{T, d}$, it follows that    
    \begin{align}
    %\sum_{i,j\in ({1\over \e}\Omega_q(2\delta)\times\cubeTk)\cap X} \e^{d-2}a_{ij}(v^{\e,q}_i-v^{\e,q}_j)^2 
    F^q_\e(u_\e^q)& = \sum_{i,j\in ({1\over \e}\Omega_q\times\cubeTk)\cap X} \e^{d-2}a_{ij}(u^{\e, q}_i-u^{\e,q}_j)^2\notag\\
    &\leq\sum_{i, j\in (\cup_{l\in{\cal I}^d_{\e,q}} Q^{l}_{T,d}\times\cubeTk)\cap X} \e^{d}a_{ij}(w^q(i)-w^q(j))^2\notag\\
    %&\leq \sum_{l\in{\cal I}^d_\e} \sum_{\e i,\e j\in (\e Q^{lK}_{K,d}\times\e\cubeTk)\cap\e X} \e^{d}a_{ij}(w_i-w_j)^2\notag\\
    &\leq \sum_{l\in{\cal I}^d_\e} \e^{d} \sum_{ i, j\in (Q^{l}_{T,d}\times\cubeTk)\cap X} a_{ij}(w^q(i)-w^q(j)^2
    %&\leq \sum_{l\in{\cal I}^d_\e}\e^d K^d(A_\hom z\cdot z + \delta)\notag\\
\leq |\Omega_q| f_\hom (z_q)+o(1),\notag
    \end{align}
as $\e\to 0$. This combined with \eqref{estOmegaq} implies that 
   \begin{equation}
   \limsup_{\e\to 0} F^q_\e(v_\e^{q,\delta})\leq |\Omega_q|f_\hom(z_q) + o(1)\label{limsupOmegaq}
   \end{equation}
       as  $\delta\to 0$.
 
\par Now, we define the recovery sequence  $v^\e$ by
    \begin{equation}
    \notag
    v^{\e}_i=v^{\e, q,\delta}_i \qquad \mbox{if } i\in \Omega_q, 
    \end{equation} 
for $q=1,\dots, \lambda$. To conclude the proof, it remains to show that, given $q_1, q_2\in \{1,\dots, \lambda\}$,  
    \begin{equation}
     \limsup_{\e\to0} \sum_{\substack{i\in({1\over \e}\Omega_{q_1}\times\cubeTk)\cap X\\j\in({1\over \e}\Omega_{q_2}\times\cubeTk)\cap X}} \e^{d-2}a_{ij}(v^{\e}_i-v^{\e}_j)^2 =o(1)
     \label{estlimsupdif}
   \end{equation}
       as  $\delta\to 0$.
 
Since $\delta>\!>\e T$, the interactions between nodes in $({1\over \e}\Omega_{q_1}(2\delta)\times\cubeTk)\cap X$ and $({1\over \e}(\Omega_{q_2}\setminus\overline{\Omega_{q_2}(\delta)})\times\cubeTk)\cap X$ or nodes in $({1\over \e}\Omega_{q_1}(2\delta)\times\cubeTk)\cap X$ and $({1\over \e}\Omega_{q_2}(2\delta)\times\cubeTk)\cap X$ do not take place. This allows to reduce    \eqref{estlimsupdif}  to the following estimate
%    \begin{equation}
%    \limsup_{\e\to0} \sum_{\substack{i\in({1\over \e}(\Omega_{q_1}\setminus\overline{\Omega_{q_1}(\delta)})\times\cubeTk)\cap X\\j\in({1\over \e}(\Omega_{q_2}\setminus\overline{\Omega_{q_2}(\delta)})\times\cubeTk)\cap X}} \e^{d-2}a_{ij}(z_{q_1}\cdot i^d-z_{q_2}\cdot\e j^d)^2 =o(1) \label{estlimsupdif1}
%   \end{equation}
    \begin{equation}
    \limsup_{\e\to0} \sum_{\substack{i\in[{1\over \e}(\Omega_{q_1}\setminus\overline{\Omega_{q_1}(\delta)})\times\cubeTk]\cap X\\j\in[{1\over \e}(\Omega_{q_2}\setminus\overline{\Omega_{q_2}(\delta)})\times\cubeTk]\cap X}} \e^{d}a_{ij}(u^\e_{i}-u^\e_{j})^2 =o(1) \label{estlimsupdif1}
   \end{equation}
as $\delta\to 0$, where we have set
$u^\e_i= u(\e i^d)$, and
used the fact that $v^{\e,q,\delta}_i=z_q\cdot\e i^d+c_d= u^\e_{i}$ if $i\in (\Omega_q\setminus\overline{\Omega_q(\delta)}\times\cubeTk)\cap X$. %Indeed, since $|z_{q_1}\cdot \e i^d-z_{q_2}\cdot\e j^d|\leq C \e|i^d-j^d|$, for some constant $C>0$, 
By the Lipschitz continuity of $u$ we deduce that 
   \begin{align}
    \sum_{\substack{i\in[{1\over \e}(\Omega_{q_1}\setminus\overline{\Omega_{q_1}(\delta)})\times \cubeTk]\cap X
    \\j\in[{1\over \e}(\Omega_{q_2}\setminus\overline{\Omega_{q_2}(\delta)})\times \cubeTk]\cap X}} \e^{d-2}a_{ij}(u^\e_{i}-u^\e_{j})^2&
    %&\quad= \sum_{i\in[{1\over \e}(\Omega_{q_1}\setminus\overline{\Omega_{q_1}(\delta)})\times \cubeTk]\cap X}
    %\sum_{j\in[{1\over \e}(\Omega_{q_2}\setminus\overline{\Omega_{q_2}(\delta)})\times \cubeTk]\cap X} \e^{d-2}a_{ij}(z_1\cdot \e i^d-z_2\cdot \e j^d)^2\notag \\
    \quad \leq C\sum_{\substack{i\in[{1\over \e}(\Omega_{q_1}\setminus\overline{\Omega_{q_1}(\delta)})\times \cubeTk]\cap X\\j\in[{1\over \e}(\Omega_{q_2}\setminus\overline{\Omega_{q_2}(\delta)})\times \cubeTk]\cap X}} \e^{d}a_{ij}| i^d-j^d|^2\notag \\
    &\quad\leq C \max\{a_{ij}\}T^k\sum_{q=1}^\lambda \Bigl|\bigcup_{l\in[{1\over \e}(\Omega_{q}\setminus\overline{\Omega_{q}(\delta)})]} \e Q^l_{d}   \Bigr| \         
\leq \ C\delta\notag        
    \end{align}
(the final  $C$ taking into account the bound for $a_{ij}$, their range, $T$ and the ${\mathcal H}^{d-1}$ measure of the union of $\partial \Omega_q$), which proves \eqref{estlimsupdif1}. 

\par Gathering estimates \eqref{limsupOmegaq} and \eqref{estlimsupdif}, we deduce that 
   \begin{equation}
   \notag
   \limsup_{\e\to 0} \sum_{i,j\in{({1\over \e}\Omega\times\cubeTk)\cap X}} \e^{d-2}a_{ij}(v^{\e}_i-v^{\e}_j)^2\leq \sum_{q=1}^{\lambda} |\Omega_q|f_\hom(z_q) + o(1).
   \end{equation}
 as $\delta\to 0$, which concludes the proof of the upper bound.

   \begin{remark}\rm Recall that the $\Gamma$-limit of a family of non-negative quadratic forms is still a non-negative quadratic form (see {\it e.g.} \cite[Theorem 11.10]{DM93}). Applying this property in our setting, we deduce that the $\Gamma$-limit $F_\hom$ of $F_\e$ is a non-negative quadratic form.  In other words, there exists a symmetric matrix $A_\hom$ such that $f_\hom(z)= A_\hom z\cdot z $, which finally gives \eqref{periodicform}.
   \end{remark}

\subsection{Convergence of minimum problems}
In this section, we deal with minimum problems with boundary data. To this end, we derive  compactness result in the case that the functionals $F_\e$ are subject to Dirichlet boundary conditions. In the discrete setting, such conditions are imposed by introducing a parameter $r\in\mathbb{N}$ and fixing the value of $u$ in a neighbourhood of the `lateral boundary' of $\Omega\times \cubeTk$, corresponding to $i^d$ in a neighbourhood of the boundary of $\Omega\subset\mathbb{R}^d$, of size $\e r$. \par For any $r>0$ and given $\varphi\in H^1 (\R^d)$, we introduce the set  
     \begin{equation*}
     {\cal C}^{\varphi, r}_{\e} (\Omega) :=\left\{u\in {\cal C}_\e(\Omega)\hspace{0.03cm}:\hspace{0.03cm} u(\e i) =\fint_{\e i^d+\e Q_d}\varphi(x^d)dx^d, \quad\mbox{if } (\e i^d + (-\e r, \e r)^d)\cap \R^d\setminus\Omega\neq\emptyset   \right\}.
     \end{equation*}
%{\color{red} Note that if $\varphi\in W^{1,\infty}(\R^d)$, in the definition above we can take $u(\e i)= \varphi(\e i^d)$}.
\par We define the functional $F^{\varphi, r}_\e$ by 
      \begin{equation}
      F^{\varphi, r}_\e (u) := F_\e(u), \quad u\in {\cal C}^{\varphi, r}_{\e} (\Omega).\notag
      \end{equation}

   \begin{thm}
   For any $\varphi\in  H^{1}(\R^d)$, let $F^\varphi$ be the functional defined by 
       \begin{equation}
       \notag
       F^\varphi(u):= \begin{dcases}
       \int_{\Omega}A_\hom \nabla u \cdot \nabla u\,dx, & u-\varphi\in H^1_0(\Omega),\\
       \infty, &\mbox{otherwise},
       \end{dcases} 
       \end{equation}
    where $A_\hom$ is given by \eqref{periodicform}. 
    Then, for any $r>0$, $F_\e^{\varphi, r}$ $\Gamma$-converge to the functional $F^\varphi$ with respect to convergence \eqref{defconvergence}.
   \end{thm}
   \begin{proof}
   We prove the $\Gamma$-liminf inequality. To that end, we prove that if $u_\e$ converges to $u$ with respect to convergence \eqref{defconvergence} and $F^{\varphi, r}_\e(u_\e)$ is equibounded, then $u-\varphi\in H^1_0(\Omega)$. First, note that if $\sup_{\e>0} F^{\varphi, r}_\e(u_\e)<\infty$, then, thanks to the coerciveness of the coefficient $a_{ij}$, we deduce that
      \begin{equation*}
      \sup_{\e\to 0} \sum_{i,j\in({1\over \e}\Omega\times\cubeTk)\cap X} \e^{d-2}(u^\e_i-u^\e_j)^2<\infty.
      \end{equation*}
   We denote by $\tilde{u}_\e$ the extension of $u_\e$ on 
    the whole $X$ defined by $\tilde{u}^\e_i = \varphi(\e i^d)$, for any $\e>0$ and outside $\Omega$. Analogously,  $\tilde{u}$ is the extension of $u$ on $\R^d$ obtained by setting $ \tilde{u}(x^d)=\varphi(x^d)$.
   Let $\Omega'$ be an open set such that $\Omega\subset\subset\Omega'$. Hence, we have that 
      \begin{align}
      \sum_{i,j\in({1\over \e}\Omega'\times Q_{T,k})\cap X} \e^{d-2}a_{ij}(\tilde{u}^\e_i-\tilde{u}^\e_j)^2 &\leq \sum_{i,j\in({1\over \e}\Omega\times\cubeTk)\cap X}\e^{d-2} a_{ij}(u^\e_i-u^\e_j)^2\notag\\
      &\quad + \sum_{i,j\in({1\over \e} (\Omega'\setminus\Omega_r)\times\cubeTk)\cap X} \e^{d-2}a_{ij}(\varphi(\e i^d)-\varphi(\e j^d))^2\leq C.\notag
      \end{align}
   Repeating similar arguments as the proof of $\Gamma$-$\liminf$ inequality of Theorem \ref{thm:homogenization} and since $\tilde{u}_\e$ converges to $\tilde{u}$, we deduce that $\tilde{u}\in H^1(\Omega')$ and hence $u-\varphi\in H^1_0(\Omega).$
   Then, invoking again Theorem \ref{thm:homogenization}, we have that
       \begin{equation*}
       \liminf_{\e\to 0} F^{\varphi, r}_\e(u_\e)=\liminf_{\e\to 0} F_\e(u_\e)\geq F^\varphi (u), 
       \end{equation*}
   as desired. \par Now, we show the $\Gamma$-limsup inequality. First, consider the case where $u\in H^1(\Omega)$ such that ${\rm supp}(u-\varphi)\subset\subset\Omega$. The general case is obtained by a density argument.\par Consider a target function $u$ such that ${\rm supp}(u-\varphi)\subset\subset\Omega$. In view of Theorem \ref{thm:homogenization}, we know that there exists a recovery sequence $u_\e$ converging to $u$ such that 
       \begin{equation}
       \notag
       \lim_{\e\to 0} F_\e(u_\e) = \int_{\Omega}A_\hom \nabla u\cdot\nabla u dx. 
       \end{equation}
   In order to modify the sequence $u_\e$ near the boundary of $\Omega$, we apply  Lemma \ref{lemma:treatmentBC} with $v_\e =u$. Hence, there exists a sequence $w_\e$ such that $w_\e$ still converges to $u$ with respect to convergence \eqref{defconvergence}, $w_\e=u$ is a neighbourhood of $\Omega$ and 
       \begin{equation}
       \limsup_{\e\to 0} F_{\e}(w_\e)\leq \limsup_{\e\to 0} F_\e(u_\e) + o(1).\notag
       \end{equation} 
   Since ${\rm supp}(u-\varphi)\subset\subset \Omega$, it follows that $w_\e$ is equal to $\varphi$ is a neighbourhood of $\Omega$, so that $F^{\varphi, r}_\e(w_\e) = F_\e (w_\e)$. We may conclude that 
       \begin{equation}
       \notag
       \limsup_{\e\to 0}F^{\varphi, r}_\e(w_\e) \leq \limsup_{\e\to 0} F_\e(u_\e) + o(1) = F^\varphi(u) + o(1),
       \end{equation}  
   which concludes the proof.
   \end{proof}
Now, we state the following result which deals with convergence of minimum problems with Dirichlet boundary data.
   \begin{prop}
   We have that
       \begin{equation}
       \notag
       \lim_{\e\to 0} \inf\{ F_\e(u)\hspace{0.03cm}:\hspace{0.03cm} u\in  {\cal C}^{\varphi, r}_\e(\Omega) \} = \min \{F_\hom (u)\hspace{0.03cm}:\hspace{0.03cm} u-\varphi \in H^1_0(\Omega) \}.
       \end{equation}
    Moreover, if $u_\e\in{\cal C}^{\varphi, r}_\e(\Omega) $ converges to $\tilde{u}$ with respect to convergence \eqref{defconvergence} and it is such that 
       \begin{equation}
       \notag
       \lim_{\e\to 0} F_{\e}(u_\e) = \lim_{\e\to 0} \inf\{ F_\e(u)\hspace{0.03cm}:\hspace{0.03cm} u\in  {\cal C}^{\varphi, r}_\e(\Omega) \},
       \end{equation}
   Then, $u$ is a minimizer for $\min \{F_\hom (u)\hspace{0.03cm}:\hspace{0.03cm} u-\varphi \in H^1_0(\Omega) \}$.
   \end{prop} 
   \begin{proof}
   We have to show  the equi-coerciveness of $F^{\varphi, r}_\e$ with respect the topology defined by \eqref{defconvergence}. To that end, consider $\{ u_\e\}\subset{\cal C}_\e^{\varphi, r}(\Omega)$ such that $\sup_{\e>0} F^{\varphi, r}_\e(u_\e)<\infty$. In view of inequality \eqref{Poincareineq} applied to $u-\varphi$, we deduce that
       \begin{align}
       \sum_{l\in\Z^d}\sum _{i,j\in (Q^l_{T,d}\times\cubeTk)\cap X}\e^d |u^\e_i-\varphi_i| &\leq C \sum_{i,j\in({1\over \e}\Omega\times\cubeTk)\cap X}\e^da_{ij}|(u^\e_i-\varphi_i)-(u^\e_j-\varphi_j)|^2\notag\\
       &\leq C F_\e(u_\e) + C\sum_{i,j\in({1\over \e}\Omega\times\cubeTk)\cap X}\e^da_{ij}|\varphi_i-\varphi_j|^2\leq C.\notag
       \end{align} 
      Hence, we may apply Proposition \ref{prop:cptresultmeanvalue} to deduce that there exists a subsequence $u_\e$ such that $\overline{u}_\e$ is converging. This concludes the proof.
      \end{proof}
The next proposition shows the Poincar\'{e} inequality for functions $u\in{\cal C}^{\varphi, r}_\e(\Omega)$. We prove it assuming that $\varphi=0$. 
      \begin{prop}
      Let $\Omega$ be a bounded open set of $\R^d$ and let $u$ be a function in ${\cal C} _\e(\Omega)$ such that $u_i=u(\e i)=0$ if $\dist(\e i^d, \partial \Omega)\leq 2\e\sqrt{d}T$. Then, there exists a constant $C>0$ such that
                \begin{equation}
                \label{Poincareineq}
                \sum_{i,j\in({1\over \e}\Omega\times\cubeTk)\cap X}\e^d|u_i|^2\leq C \sum_{i,j\in({1\over \e}\Omega\times\cubeTk)\cap X}\e^da_{ij}|u_i-u_j|^2,
                \end{equation}
      where $C$ is of order of $[{\rm diam}(\Omega)]^2$.
      \end{prop}
      
      \begin{proof}
      We identify $u$ with its extension to $(\R^d\times\cubeTk)\cap X$ which is equal to $0$ outside $(\Omega\times\cubeTk)\cap X$.  Due to the boundedness of $\Omega$, there exists  $M>0$ such that $\Omega\subset [0,M)^d$ and $({1\over \e}[0, M)^d\times\cubeTk)\cap X$ contains a path joining two arbitrary nodes $i,j\in ({1\over \e}\Omega\times\cubeTk)\cap X$. \par Fix $i\in ({1\over \e}\Omega\times\cubeTk)\cap X$ and let $j$ be a node such that $\dist (\e j^d, \partial\Omega)\leq 2\e\sqrt{d}T$ and $i^d-j^d=\lambda T e_1$, where, without lost of generality, we may assume that  $\lambda$ is a positive integer. Note that $\lambda$ depends on the fixed node $i$ and it is of order $MT^{-1}\e^{-1}$.      
      Let  $l_i$ and $l_j$ be two indices in $\Z^d$ such that $i\in Q^{l_i}_{T,d}$ and $j\in Q^{l_j}_{T,d}$. Let $S_{T,d}^\lambda$ be the union of $(\lambda+1)$ neighbour cubes joining $Q^{l_i}_{T,d}$ and $Q^{l_j}_{T,d}$ such that each two consecutive cubes having one face in common. In other words,
         \begin{equation}
         \notag
         S_{T,d}^\lambda := \bigcup_{q=0}^{\lambda} Q^{l_q}_{T,d},
         \end{equation}
      where $l_q=l_{j}+qTe_1$, for $q=1,\dots, \lambda$ and $l_0=l_j$. Since $X$ is connected, there exists a path of nodes $\{j_q\}_{q=0}^{\lambda}$ joining $j_0=j$ and $j_\lambda=i$  such that it is contained in $S_{T,d}^\lambda+(-T, T)^d$, $j_q\in(Q^{l_q}_{T,d}\times\cubeTk)\cap X$ and $(j_q,j_{q+1})\in{\cal E}$. Such a path can be built repeating periodically the path joining $j\in (Q^{l_0}_{T,d}\times\cubeTk)\cap X$ and $j+Te_1\in(Q^{l_1}_{T,d}\times\cubeTk)\cap X $. Since $u_j=0$, we have that 
              \begin{equation}
              u_i =\sum_{q=1}^{\lambda}(u_{j_q}-u_{j_{q-1}}).\notag 
              \end{equation}
      Hence, an application of the Jensen inequality leads us to  
                \begin{equation}
                |u_i|^2\leq\lambda \sum_{q=1}^{\lambda} |u_{j_q}-u_{j_{q-1}}|^2.\notag
                \end{equation}
            Summing over $ i\in({1\over \e}\Omega\times\cubeTk)\cap X$, we get 
                \begin{align}
                \sum_{i\in ({1\over \e}\Omega\times\cubeTk)\cap X}\e^d |u_i|^2& %=\sum_{q=1}^{\lambda}\sum_{i\in (Q^{l_q T}_{T,d}\times\cubeTk) \cap X}\e^d |u_i|^2 \notag\\
                \leq \lambda \sum_{q=1}^{\lambda}\sum_{i,j\in ({1\over \e}\Omega\times\cubeTk)\cap X}\e^d|u_{j_q}-u_{j_{q-1}}|^2\notag\\
                &\leq C\lambda^2 \sum_{i,j\in ({1\over \e}\Omega\times\cubeTk)\cap X}\e^da_{ij}|u_{j_q}-u_{j_{q-1}}|^2\notag,
                %\lambda \sum_{q=1}^{\lambda}\sum_{i,j\in S^\lambda_{T,d}+(-\delta,\delta)^d\times\cubeTk\cap X}\e^da_{ij}|u_{j_q}-u_{j_{q-1}}|^2\notag
                \end{align}
where the constant $C$ takes into account the fact that the possible multiplicity of the paths containing the connection joining $j_{q-1}$ and $j_q$, which is anyhow uniformly bounded. Recalling that $\lambda$ is of order $MT^{-1}\e^{-1}$, we get the inequality \eqref{Poincareineq}, as desired.
\begin{comment}
            Since $S^\lambda_{T,d}+(-\delta, \delta)^d$ contained the path of nodes $\{j_q\}_{q=0}^{\lambda}$ and due to the property of the coefficient $a_{ij}$, we have that 
                \begin{align}
                \sum_{i,j\in (S^\lambda_{T,d}+(-\delta, \delta)^d\times\cubeTk)\cap X}\e^d|u_{j_q}-u_{j_{q-1}}|^2&\leq C\sum_{i,j\in (S^\lambda_{T,d}+(-\delta, \delta)^d\times\cubeTk)\cap X}\e^da_{ij}|u_{j_q}-u_{j_{q-1}}|^2,\notag
                \end{align}
            which implies that 
               \begin{align}
               \sum_{q=1}^{\lambda}\sum_{i\in (Q^{l_q T}_{T,d}\times\cubeTk) \cap X}\e^d |u_i|^2&\leq C\lambda\sum_{q=1}^{\lambda}\sum_{i,j\in (S^\lambda_{T,d}+(-\delta, \delta)^d\times\cubeTk)\cap X}\e^da_{ij}|u_{j_q}-u_{j_{q-1}}|^2.\notag
               \end{align}    
             We sum over $l\in\Z^d$ obtaining that
               \begin{align}
               \sum_{l\in\Z^d}\sum_{i\in (Q^{l T}_{T,d}\times\cubeTk) \cap X}  \e^d|u_i|^2  
               &\leq C \lambda^2\sum_{i,j\in ({1\over \e}\Omega\times\cubeTk)\cap X} \e^{d}a_{ij}|u_i-u_j|^2\notag\\
               &\leq C {M^2\over T^2}\sum_{i,j\in ({1\over \e}\Omega\times\cubeTk)\cap X} \e^{d-2}a_{ij}|u_i-u_j|^2,\notag
               \end{align}
          where in the last inequality we have used the fact that $\lambda$ is of order $MT^{-1}\e^{-1}$.
\end{comment}          
      \end{proof}

\section{Examples}\label{Examples}
In this section, we exhibit some examples of the possible geometries of the set $X$. We also  compute  the homogenized matrix $A_\hom$ given by formula \eqref{periodicform}.  In the examples below,  we think of $X$ as a subset of  $\Z^{d+k}$ where $d=1$ is identified with the horizontal direction and $k=1$ or $2$. 
Since $d=1$ the homogenized matrix actually reduces to a single coefficient giving the homogenized energy density $A_\hom z^2$. 
%Note that since $A_\hom$ turns out to be  positively homogeneous of degree $2$, it is sufficient to compute 
  %   \begin{align}
  %   A_\hom  ={1\over T^d}\inf\biggl\{\sum_{i\in(\cubeTd\times\cubeTk)\cap X}\sum_{j\in( \R^d\times\cubeTk)\cap X}a_{ij}(u_i-u_j)^2\hspace{0.03cm}:\hspace{0.03cm}u_i- z\cdot i^d\quad  \mbox{is } T\mbox{-periodic in }e_1,\dots,e_d \biggr\}. \label{homgenAhom}
  %   \end{align}

In all the following examples the value of the non-zero coefficients $a_{ij}$ is always $1$, and the corresponding connections are represented by solid lines in the figures.

      \begin{figure}[h]
        \centering
        \begin{minipage}{0.35\textwidth}
        \centering
        \includegraphics[scale=0.90]{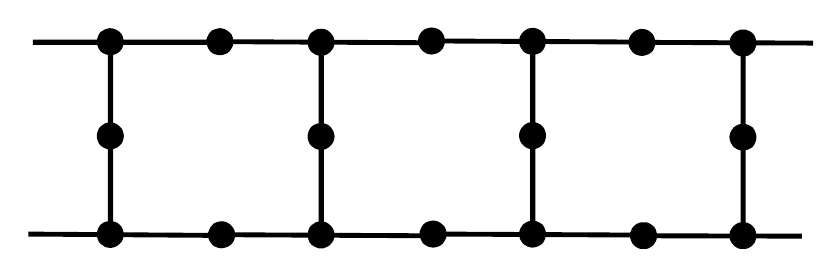}
        \caption*{(a)}
        \end{minipage} %
        \hspace{2.3cm}
        \begin{minipage}{0.35\textwidth}
        \centering
        \includegraphics[scale=0.90]{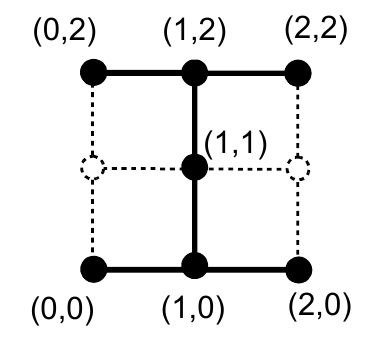}
        \caption*{(b)}
        \end{minipage} %
         \caption{\small Figure (a) shows $X$ and Figure (b) shows the periodicity cell $(Q_{2,d}\times Q_{2,k})\cap X$.} 
         \label{fig:example1}
        \end{figure}

\begin{example}
  \label{Ex1}
\rm  Let $X$ be the set pictured in Figure \ref{fig:example1}(a). Here, we have that $d=k=1$ and the period $T$ is equal to $2$. Figure \ref{fig:example1}(b) shows a periodicity cell.
The geometry of the set $X$ can be thought as the discrete version of a perforated domain. Indeed, note that  nodes $(0,1)$ and $(2,1)$ in Figure \ref{fig:example1}(b) are missing.
%  \par To compute $A_\hom$, first note that in view of homogeneity, we may assume that $z=1$ in the formula \eqref{periodicform}. We may also assume  minimizer $\tilde{u}$ of $A_\hom$ be such that $\tilde{u}(0,0)=0$. The periodicity condition implies that $\tilde{u}(2,0)= 2$ and $\tilde{u}(2,2)= \tilde{u}(0,2)+2$.
%  Set $\tilde{u}(1,0)=u^-$ and $\tilde{u}(1,2)=u^+$. An application of the Jensen inequality to the energy leads us to 
%        \begin{align}
%       (u^-)^2&+(2-u^{-})^2 + (u^- -u(1,1))^2 + (u(1,1)-u^+)^2+( u^+-\tilde{u}(0,2))^2 + (u^+-\tilde{u}(0,2)-2 )^2\notag\\
%        &\geq {1\over 2}(2-2u^-)^2+{1\over 2}(u^- + u^+-2u(1,1))^2 + {1\over 2}(2u^+-\tilde{u}(0,2)-2 )^2.\notag
%        \end{align} 
%  This suggests that the minimizer $\tilde{u}$ is such that $u^-= 1$, $u^+=\tilde{u}(0,2)+1$ and  $\tilde{u}(1,1)= 1+{\tilde{u}(0,2)\over 2}$. Hence, minimizing the energy
%      $
%        4+ {\tilde{u}(0,2)^2\over 2},\notag
%       $
%   we find that $\tilde{u}(0,2)=0$. 
A minimizer $\tilde{u}$  for \eqref{periodicform} is given by $\tilde{u}(0,0)=\tilde{u}(0,2)=0$, $\tilde{u}(1,0)=\tilde{u}(1,1)=\tilde{u}(1,2)=z$ and $\tilde{u}(2,0)=\tilde{u}(2,2)=2z$, so that $A_\hom = 4$.  
  \end{example}
  
  In the next four examples $d=k=1$, the set $X$ is always simply $\mathbb Z\times \{0,1\}$ and the period $T$ is $1$, but the set $\mathcal E$ is such that the graph cannot be directly seen as a discretization of a thin film in the continuum parameterized as a subgraph of a function of one real variable .

 \begin{figure}[h]
  \centering
  \begin{minipage}{0.4\textwidth}
  \centering
  \includegraphics[scale=0.40]{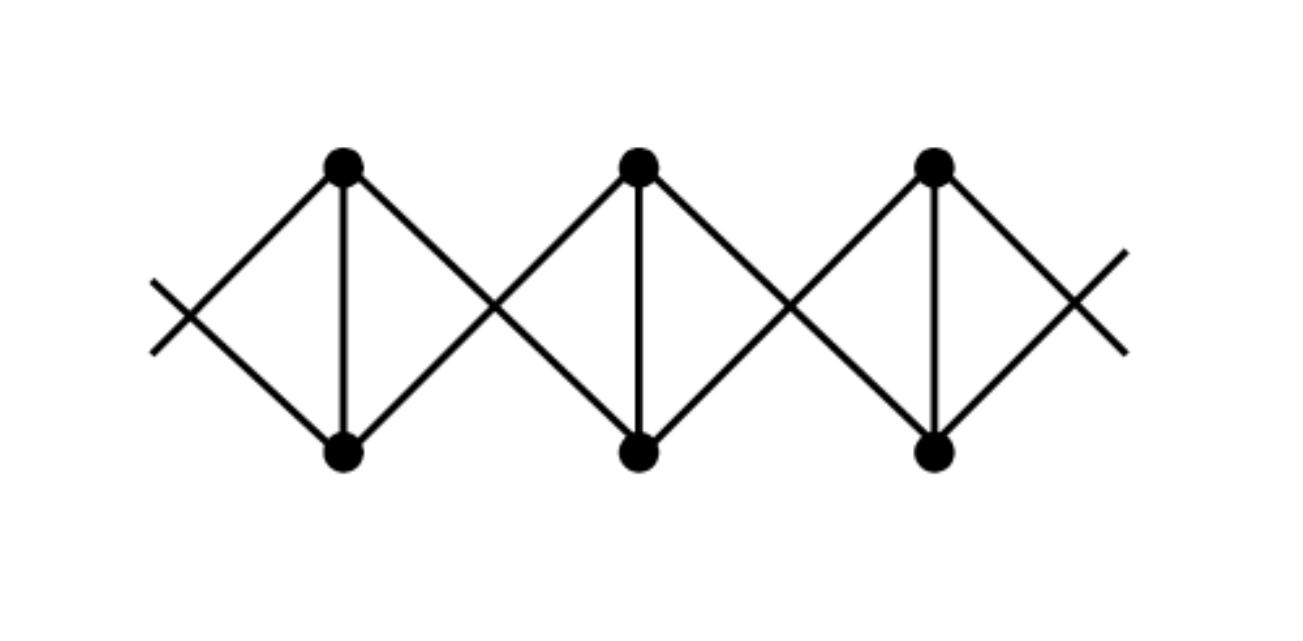}
  \caption*{(a)}
  \end{minipage} %
  \hspace{1cm}
  \begin{minipage}{0.4\textwidth}
  \centering
  \includegraphics[scale=0.40]{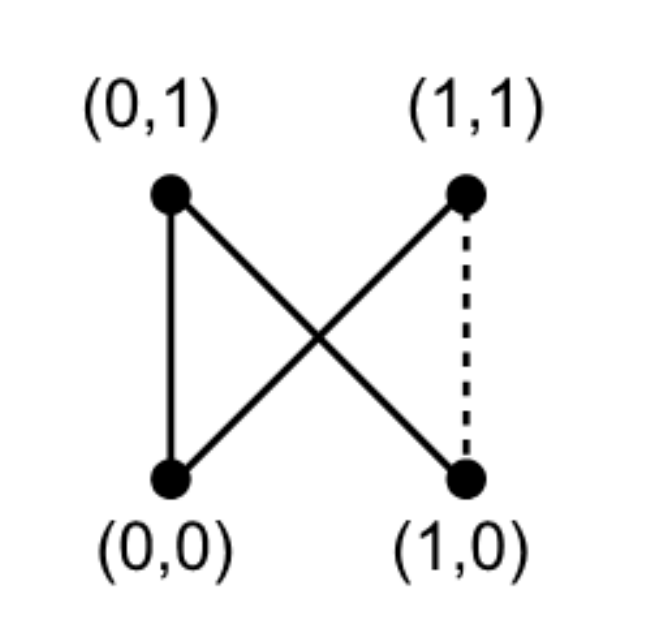}
  \caption*{(b)}
  \end{minipage} %
   \caption{\small Figure (a) shows $X$ and Figure (b) shows the periodicity cell $(Q_{1,d}\times Q_{1,k})\cap X$.} 
   \label{fig:example2}
  \end{figure}

  \begin{example}
  \label{Ex2}
  \rm Let $X$ be as drawn in Figure \ref{fig:example2}(a). 
In this case $\mathcal E$  contains all `cross-connections' between points of $X$.
%  {\rm {\color{red} The set $X$ can be thought as a structure of conducting rods intersecting in the centrer of the square $[0,1)\times [0,1)$. Note that $X$ is the entire lattice $\Z^2$ where some links among nodes are been removed. However, the structure of $X$ remains connected and hence $X$ satisfies our assumptions.}
% In order to compute the homogenized matrix $A_\hom$, we repeat similar arguments as in Example \ref{Ex1}, obtaining that  and 
 The minimizer $\tilde{u}$ for $A_\hom z^2$ is  $\tilde{u}(0,0)=\tilde{u}(1,1)=0$ and $\tilde{u}(1,0)= \tilde{u}(0,1)=z$, so that $A_\hom = 4$.
  %As in the previous example, we may assume that the minimizer $\tilde{u}$ satisfying $\tilde{u}(0,0) = 0$. By  the periodicity condition, we also know that $\tilde{u}(0,1)=1$. Set $u\tilde{u}(1,0) = w$ which implies  $\tilde{u}(1,1)= w+1$.  Hence, minimizing the energy $2w^2 + 4(w-1)^2+2(w-2)^2$, we get that  $\tilde{u}$ is such that $\tilde{u}(0,0)=\tilde{u}(1,1)=0$ and $\tilde{u}(1,0)= \tilde{u}(0,1)=z$. We conclude that $A_\hom = 4$.}
  \end{example}
 
     \begin{figure}[h]
    \centering
    \begin{minipage}{0.4\textwidth}
    \centering
    \includegraphics[scale=0.40]{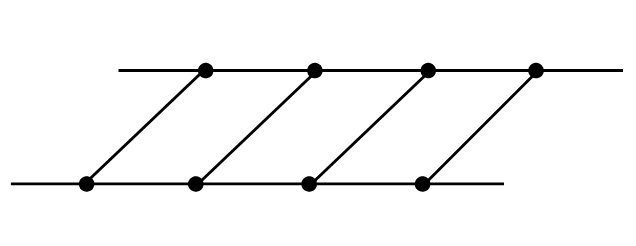}
    \caption*{(a)}
    \end{minipage} %
    \hspace{1cm}
    \begin{minipage}{0.4\textwidth}
    \centering
    \includegraphics[scale=0.40]{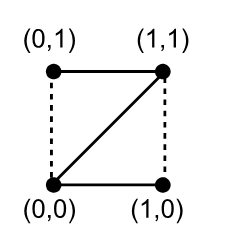}
    \caption*{(b)}
    \end{minipage} %
     \caption{\small Figure (a) shows $X$ and Figure (b) shows the periodicity cell $(Q_{1,d}\times Q_{1,k})\cap X$.} 
     \label{fig:example3}
    \end{figure}
   
   \begin{example}
 \rm Consider $X$ as drawn in Figure \ref{fig:example3}(a).  
 Here, the graph is analog to a nearest-neighbour thin film, but with a translation of a unit of one of the two copies of $\mathbb Z$, which again makes this geometry not immediately seen as a discretization of a continuum thin film.
%  
%  
%  {\color{red} The structure of $X$ appears skew with respect to the ones analysed in Examples \ref{Ex2} and \ref{Ex3}.  Note that the periodicity cell (see Figure \eqref{fig:example3} (b)) is similar to the one in Example \ref{Ex3} where the link connecting $(0,0)$ and $(0,1)$ is removed.}{\rm  With the same argument of Example \ref{Ex1}, 
 The $1$-periodic minimizer $\tilde{u}$ for $A_\hom z^2$ is given by $\tilde{u}(0,0)=\tilde{u}(1,1)=0$, $\tilde{u}(1,0)=z$ and $\tilde{u}(0,1)=-z$ and the homogenized coefficient is $A_\hom=4$.
   % such that $\tilde{u}(0,0) =0$ which implies that $\tilde{u}(1,0)=1$. Set $\tilde{u}(1,1)=w$. Hence,  $\tilde{u}(0,1)= w-1$. Minimizing the energy $4+4w^2$, we obtain that the minimizer $\tilde{u}$ is given by $\tilde{u}(0,0)=\tilde{u}(1,1)=0$, $\tilde{u}(1,0)=z$ and $\tilde{u}(0,1)=-z$. The homogenized matrix $A_\hom$ is equal to $4$.}
  \end{example}

    \begin{figure}[h]
    \centering
    \begin{minipage}{0.4\textwidth}
    \centering
    \includegraphics[scale=0.40]{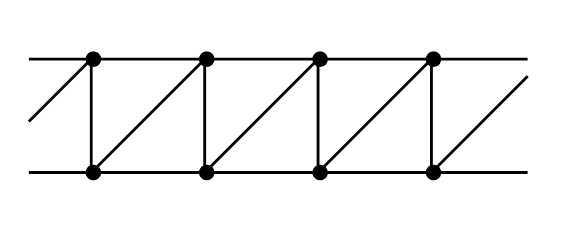}
    \caption*{(a)}
    \end{minipage} %
    \hspace{1cm}
    \begin{minipage}{0.4\textwidth}
    \centering
    \includegraphics[scale=0.40]{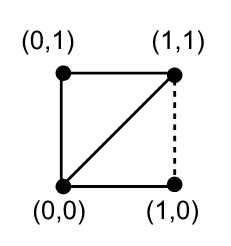}
    \caption*{(b)}
    \end{minipage} %
     \caption{\small Figure (a) shows $X$ and Figure (b) shows the periodicity cell $(Q_{1,d}\times Q_{1,k})\cap X$.} 
     \label{fig:example6}
    \end{figure}
  
  \begin{example}
  \label{Ex3}\rm
Consider $X$ as drawn in Figure \ref{fig:example6}(a). Here the set of connections has the structure of a triangular lattice.
The minimizer $\tilde{u}$ for $A_\hom z^2$ is given by $\tilde{u}(0,0)=0$, $\tilde{u}(1,0)=z$, $\tilde{u}(0,1)=-1/2z$ and $\tilde{u}(1,1)=1/2z$. The homogenized coefficient is $A_\hom=5/2$. 
   \end{example}

       \begin{figure}[h]
           \centering
           \begin{minipage}{0.4\textwidth}
           \centering
           \includegraphics[scale=0.26]{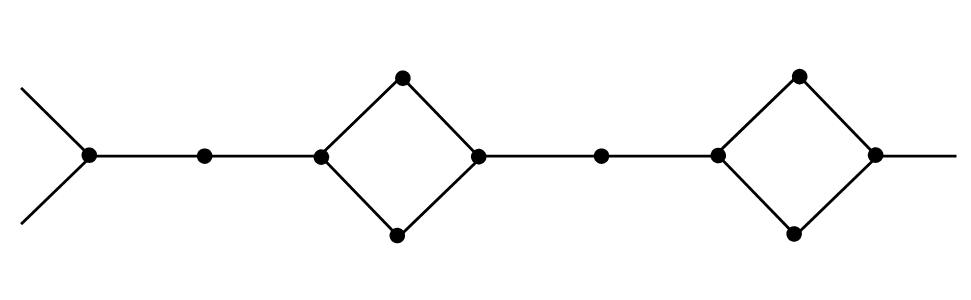}
           \caption*{(a)}
           \end{minipage} %
           \hspace{1cm}
           \begin{minipage}{0.4\textwidth}
           \centering
           \includegraphics[scale=0.40]{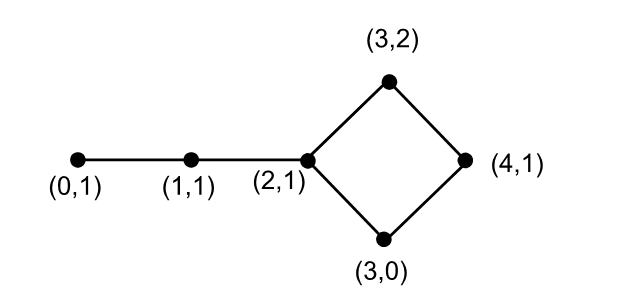}
           \caption*{(b)}
           \end{minipage} %
            \caption{\small Figure (a) shows $X$ and Figure (b) shows the periodicity cell $(Q_{4,d}\times Q_{4,k})\cap X$.} 
            \label{fig:example4}
           \end{figure}

  \begin{example}
    \rm Let $X$ be the set pictured in Figure \ref{fig:example4}(a), where $d=k=1$  and the period $T$ is equal to $4$. The set $X$ is a subset of $\Z\times \{0,1,2\}$. Such a set $X$ can be though as a discrete layered media, whose conductivity is equal to $1$ along the straight lines, while in the part corresponding to the rhombus structure the effective conductivity is $2$.
           
%    {\rm  Now, we compute the homogenized matrix $A_\hom$ and the corresponding minimizer $\tilde{u}$. We assume that $\tilde{u}$ is such that $\tilde{u}(0,1)=0$ which implies $\tilde{u}(4,1)=4$. Set $\tilde{u}(2,1)=2w$. By Jensen's inequality, we have that
%       \begin{align}
%       {1\over 2}&\left[(\tilde{u}(1,1))^2 + (\tilde{u}(1,1)-2w)^2+(\tilde{u}(3,2)-2w)^2+2(\tilde{u}(3,0)-2w)^2+2(\tilde{u}(3,2)-4)^2+2(\tilde{u}(3,0)-4)^2\right]\notag\\
%       &\geq {1\over 4} (2\tilde{u}(1,1)-2w)^2+{1\over 4}(2w+4-2\tilde{u}(3,2))^2+{1\over 4}(2w+4-2\tilde{u}(3,0))^2.\notag
%       \end{align}
%    This suggest that  $\tilde{u}$ satisfies the conditions $\tilde{u}(1,1)=w$ and  $\tilde{u}(3,0)=\tilde{u}(3,2)=2+w$.  Now, we minimize the energy $w^2+2(w-2)^2$ and we conclude that 
The minimizer $\tilde{u}$ for $A_\hom z^2$ is given by $\tilde{u}(0,0)=z$, $\tilde{u}(1,1)=4z/3$, $\tilde{u}(2,1)=8z/3$, $\tilde{u}(3,0)=\tilde{u}(3,2)= 10z/3$ and $\tilde{u}(4,1)=4z$ and $A_\hom= 8/3$.
    \end{example}
    
           \begin{figure}[h]
          \centering
          \begin{minipage}{0.4\textwidth}
          \centering
          \includegraphics[scale=0.30]{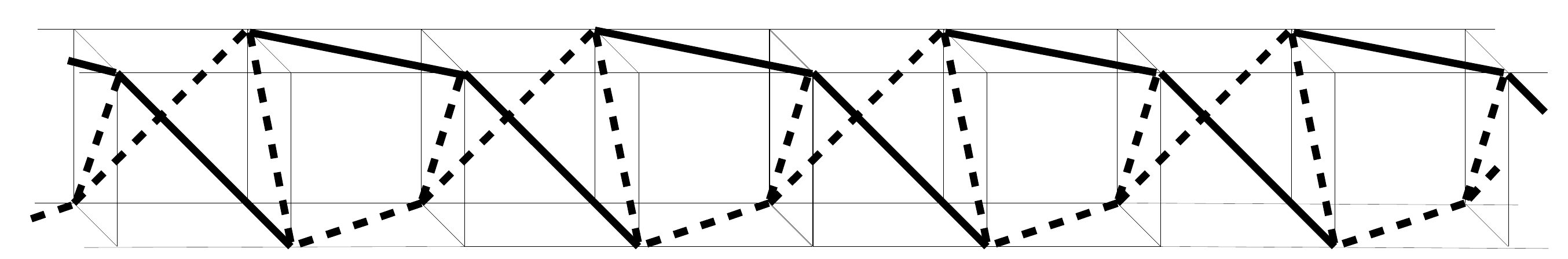}
          \caption*{(a)}
          \end{minipage} %
          \hspace{2.4cm}
          \begin{minipage}{0.4\textwidth}
          \centering
          \includegraphics[scale=0.30]{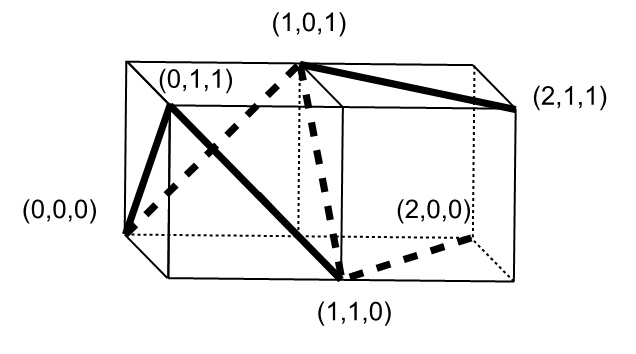}
          \caption*{(b)}
          \end{minipage} %
           \caption{\small Figure (a) shows $X$ and Figure (b) shows the periodicity cell $(Q_{2,d}\times Q_{2,k})\cap X$.} 
           \label{fig:example5}
          \end{figure}
    \begin{example}\rm
We consider the set $X$ drawn in Figure \ref{Fig5larger_vertical}. To uniform the notation introduced in this section, we rotate $X$, obtaining the structure pictured in Figure \ref{fig:example5}(a). Here $d=1$ and $k=2$. The period $T$ is equal to $2$ and the periodicity cell is drawn in Figure \ref{fig:example6}(b). 
%We compute $A_\hom$ and the corresponding minimizer $\tilde{u}$.  Assume that $\tilde{u}$ be $\tilde{u}(0,0,0)=0$. By the periodicity condition, we deduce that $\tilde{u}(2,0,0)=2$. Set $\tilde{u}(1,1,0)=v$ and $\tilde{u}(1,0,1)=w$. Thanks to the Jensen inequality, we have that 
%      \begin{align}
%      &(\tilde{u}(0,1,1))^2+w^2+(\tilde{u}(0,1,1)-v)^2 +(v-w)^2+(v-2)^2+(\tilde{u}(0,1,1)+2-w)^2\notag\\
%      &\quad\geq (\tilde{u}(0,1,1))^2 + {1\over 2}(\tilde{u}(0,1,1)+2-2w)^2 +{1\over 2}(\tilde{u}(0,1,1) +2-2v)^2 + (v-w)^2.\notag
%      \end{align}
%     This implies that  $\tilde{u}$ is such that $v=w={\tilde{u}(0,1,1)\over 2} +1$. Then, we minimize the  energy given by
%         \begin{align}
%         (\tilde{u}(0,1,1))^2 +2 \left({\tilde{u}(0,1,1)\over 2} +1 \right)^2+2\left({\tilde{u}(0,1,1)\over 2} -1 \right)^2\notag
%         \end{align}
%      obtaining that $\tilde{u}(0,1,1)=0$. Hence, the minimizer $\tilde{u}$ is such that $\tilde{u}(0,0,0)= \tilde{u}(0,1,1)=0$, $\tilde{u}(1,0,1)=\tilde{u}(1,1,0)=z$ and $\tilde{u}(2,0,0)=\tilde{u}(2,1,1)=2z$. We conclude that  
The structure is actually the same as that in Example \ref{Ex2} but transposed to a three-dimensional setting, and $A_\hom= 4$. Note that in this case the solid lines representing the connections do not intersect and they have all the same length, so that they can also be interpreted as a system of homogeneous conducting rods. 

    \end{example}
  
  \subsection*{Acknowledgments.}
The authors acknowledge the MIUR Excellence Department Project awarded to the Department of Mathematics, University of Rome Tor Vergata, CUP E83C18000100006.

\bibliographystyle{plain}

\end{document}